\documentclass[11pt]{article}
\usepackage{amsfonts,a4wide}
\usepackage{amssymb}
\usepackage{graphicx}
\usepackage{mathtools}
\usepackage{enumerate}
\newtheorem{theorem}{Theorem}
 \newtheorem{corollary}[theorem]{Corollary}
 \newtheorem{lemma}[theorem]{Lemma}
 \newtheorem{proposition}[theorem]{Proposition}

 \newtheorem{definition}[theorem]{Definition}
\newtheorem{remark}[theorem]{Remark}
\newtheorem{example}[theorem]{Example}

\newcommand{\diag}{\mbox{\rm diag}}

\newcommand{\rank}{\mbox{\rm rank}}

\newcommand {\proof} {\par{\it Proof}. \ignorespaces}
\newcommand {\eproof}
      {\space
        {\ \vbox{\hrule\hbox{\vrule height1.3ex\hskip0.8ex\vrule}\hrule}}
        \par}
  \DeclareMathOperator{\ck}{ker}

\DeclareMathOperator{\sing}{sing}
\DeclareMathOperator{\Str}{\mathcal P}
\DeclareMathOperator{\hi}{hi}
\DeclareMathOperator{\inst}{inst}

\newcommand{\Real}{\mathbb{R}}
\newcommand{\Comp}{\mathbb{C}}

\newcommand{\eps}{\varepsilon}

\newcommand{\set}[1]{\left\{#1\right\}}

\newcommand{\norm}[1]{\left\Vert#1\right\Vert}

\catcode`@=11
\font\tenex=cmex10 
\newdimen\p@renwd
\setbox0=\hbox{\tenex B} \p@renwd=\wd0 
\def\bmat#1{\begingroup \m@th
  \setbox\z@\vbox{\def\cr{\crcr\noalign{\kern2\p@\global\let\cr\endline}}%
    \ialign{$##$\hfil\kern2\p@\kern\p@renwd&\thinspace\hfil$##$\hfil
      &&\quad\hfil$##$\hfil\crcr
      \omit\strut\hfil\crcr\noalign{\kern-\baselineskip}%
      #1\crcr\omit\strut\cr}}%
  \setbox\tw@\vbox{\unvcopy\z@\global\setbox\@ne\lastbox}%
  \setbox\tw@\hbox{\unhbox\@ne\unskip\global\setbox\@ne\lastbox}%
  \setbox\tw@\hbox{$\kern\wd\@ne\kern-\p@renwd\left[\kern-\wd\@ne
    \global\setbox\@ne\vbox{\box\@ne\kern2\p@}%
    \vcenter{\kern-\ht\@ne\unvbox\z@\kern-\baselineskip}\,\right]$}%
  \null\;\vbox{\kern\ht\@ne\box\tw@}\endgroup}
\catcode`@=12
\providecommand{\norm}[1]{\left\Vert#1\right\Vert}
\newcommand {\mycomment}[1]{} 
\newcommand {\mat}  [1] {\left[\begin{array}{#1}}
\newcommand {\rix}      {\end{array}\right]}

\newcommand {\ve}{\varepsilon}
\catcode`@=11     
\font\tenex=cmex10 
\newdimen\p@renwd
\setbox0=\hbox{\tenex B} \p@renwd=\wd0 
\def\bmat#1{\begingroup \m@th
  \setbox\z@\vbox{\def\cr{\crcr\noalign{\kern2\p@\global\let\cr\endline}}%
    \ialign{$##$\hfil\kern2\p@\kern\p@renwd&\thinspace\hfil$##$\hfil
      &&\quad\hfil$##$\hfil\crcr
      \omit\strut\hfil\crcr\noalign{\kern-\baselineskip}%
      #1\crcr\omit\strut\cr}}%
  \setbox\tw@\vbox{\unvcopy\z@\global\setbox\@ne\lastbox}%
  \setbox\tw@\hbox{\unhbox\@ne\unskip\global\setbox\@ne\lastbox}%
  \setbox\tw@\hbox{$\kern\wd\@ne\kern-\p@renwd\left[\kern-\wd\@ne
    \global\setbox\@ne\vbox{\box\@ne\kern2\p@}%
    \vcenter{\kern-\ht\@ne\unvbox\z@\kern-\baselineskip}\,\right]$}%
  \null\;\vbox{\kern\ht\@ne\box\tw@}\endgroup}
\catcode`@=12    
\newif\ifMDlatex
\catcode`@=11     
\def\MD@us#1{\csname#1style\endcsname}
\def\MD@uf#1{\csname#1font\endcsname}
\def\MD@t{text}
\def\MD@s{script}
\def\MD@ss{scriptscript}
\newdimen\MD@unit
\MD@unit\p@
\def\MD@changestyle#1{
  \relax\MD@unit0.1\fontdimen6\MD@uf{#1}0
  \everymath\expandafter{\the\everymath\MD@us{#1}}
}
\def\MD@dot{$\m@th\ldotp$}
\def\MD@palette#1{\mathchoice{#1\MD@t}{#1\MD@t}{#1\MD@s}{#1\MD@ss}}
\def\MD@ddots#1{{\MD@changestyle{#1}%
  \mkern1mu\raise7\MD@unit\vbox{\kern7\MD@unit\hbox{\MD@dot}}%
  \mkern2mu\raise4\MD@unit\hbox{\MD@dot}%
  \mkern2mu\raise \MD@unit\hbox{\MD@dot}\mkern1mu}}%
\def\MD@iddots#1{{\MD@changestyle{#1}%
  \mkern1mu\raise \MD@unit\hbox{\MD@dot}%
  \mkern2mu\raise4\MD@unit\hbox{\MD@dot}%
  \mkern2mu\raise7\MD@unit\vbox{\kern7\MD@unit\hbox{\MD@dot}}}}%
\def\MD@vdots#1{\vbox{\MD@changestyle{#1}%
    \baselineskip4\MD@unit\lineskiplimit\z@
    \kern6\MD@unit\hbox{\MD@dot}\hbox{\MD@dot}\hbox{\MD@dot}}}%
\ifMDlatex
  \DeclareRobustCommand\ddots{\mathinner{\MD@palette\MD@ddots}}%
  \DeclareRobustCommand\iddots{\mathinner{\MD@palette\MD@iddots}}%
  \DeclareRobustCommand\vdots{\mathinner{\MD@palette\MD@vdots}}%
\else
  \def\ddots{\mathinner{\MD@palette\MD@ddots}}%
  \def\iddots{\mathinner{\MD@palette\MD@iddots}}%
  \def\vdots{\mathinner{\MD@palette\MD@vdots}}%
\fi
\catcode`@=12    

\newcommand{\matp}[1]{\begin{bmatrix} #1 \end{bmatrix}}

\usepackage{color}

\begin{document}
\title
{Distance problems for dissipative Hamiltonian systems and related matrix polynomials}
\author{C. Mehl \footnotemark[3]~\footnotemark[1]
\and V. Mehrmann\footnotemark[3]~\footnotemark[1]
\and  M. Wojtylak \footnotemark[2]~\footnotemark[4]
}
\maketitle

\begin{center} \emph{Dedicated to Paul Van Dooren on the occasion of his 70th birthday} \end{center}

\begin{abstract}
We study the characterization of several distance problems for linear differential-algebraic systems with dissipative Hamiltonian
structure. Since all models are only approximations of reality and data are always inaccurate, it is an important question whether a given
model is close to a 'bad' model that could be considered as ill-posed or singular. This is usually done by computing a distance to the
nearest model with such properties. We will discuss the distance to  singularity and the
distance to the nearest high index problem for dissipative Hamiltonian systems. While for general unstructured differential-algebraic
systems the characterization of these distances are partially open problems, we will show that for dissipative Hamiltonian systems and
related matrix polynomials there exist explicit characterizations that can be implemented numerically.
\end{abstract}


\noindent

\renewcommand{\thefootnote}{\fnsymbol{footnote}}
\footnotetext[3]{
Institut f\"ur Mathematik, MA 4-5, TU Berlin, Str. des 17. Juni 136,
D-10623 Berlin, FRG.
\texttt{$\{$mehl,mehrmann$\}$@math.tu-berlin.de}.
}

\footnotetext[2]{Instytut Matematyki, Wydzia\l{} Matematyki i Informatyki,
Uniwersytet Jagiello\'nski, Krak\'ow, ul. \L ojasiewicza 6, 30-348 Krak\'ow, Poland
   \texttt{michal.wojtylak@uj.edu.pl}.}
\footnotetext[4]{   Supported by the Alexander von Humboldt Foundation.}
\footnotetext[1]{
Partially supported by {\it Deutsche Forschungsgemeinschaft} through
the  Excellence Cluster {\sc Math$^+$} in Berlin and Priority Program 1984 'Hybride und multimodale Energiesysteme:
 Systemtheoretische Methoden für die Transformation und den Betrieb komplexer Netze'.
}
\renewcommand{\thefootnote}{\arabic{footnote}}

\section{Introduction}
We study several distance problems for linear systems of \emph{differential-algebraic equations} (DAEs)
of the form
\begin{equation}\label{dH}
E\dot x = (J-R)Qx,
\end{equation}
with constant coefficient matrices  $E,Q,J, R\in\Real^{n,n}$ and a differentiable state function $x:\Real \to \Real^n$, see
also \cite{BeaMXZ18,GilMS18,JacZ12,MehMW18,MehM19,Sch13,SchM02,SchM18} for definitions and a detailed analysis of such systems in
different generality and their relation to the more general \emph{port-Hamiltonian systems}. A system of the above form  is called
linear time-invariant \emph{dissipative Hamiltonian (dH) differential-algebraic equation 
(dHDAE) system} if
\begin{equation}\label{mm5-constr}
E^\top Q \geq 0,\quad J=-J^\top, \quad  R= R^\top\geq 0,
\end{equation}
where $A^\top$ denotes the transpose of a matrix $A$ and for a symmetric matrix $A=A^\top$ by $A>0$ ($A\geq 0$) we denote that $A$
is positive definite (positive semidefinite).  Such dHDAE systems
 generalize linear time-invariant ordinary \emph{dissipative Hamiltonian systems} (the case where $E=I$) is the identity) and linear
 \emph{Hamiltonian systems}, the case that $E=I$ and  $R=0$. The associated quadratic \emph{Hamiltonian} is given by
 $\mathcal H(x)= \frac 12 x^\top E^\top Qx $ and satisfies  the \emph{dissipation inequality}
 $ {\mathcal H}\big(x(t_1)\big)-{\mathcal H}\big(x(t_0)\big) \leq 0$ for $t_1\geq t_0$.
 In many applications the matrix $Q$ can be chosen to be the identity
matrix, i.e., $Q=I$, see \cite{BeaMXZ18,GilMS18,MehM19,Sch13,SchJ14}, and this is the case that we study in this paper. In
Section~\ref{sec:noDQ} we provide an analysis how the general case \eqref{dH} can be transformed to this situation.

The system properties of \eqref{dH} (with $Q=I$) can be analyzed by investigating the corresponding dH matrix pencil
\begin{equation}\label{dHwoQ}
L(\lambda):=\lambda E-(J-R).
\end{equation}
In our analysis we focus on systems with real coefficients. Some of our results can also easily be extended to
the case of complex
coefficients, but in some occasions we make explicit or implicit use of the fact that skew-symmetric matrices have a zero
diagonal which is not true for skew-Hermitian matrices.

In many practical cases, see, e.g., \cite{BeaMXZ18,GramQSW16}, the underlying system  is of second order form
$M\ddot x-(G-D)\dot x+Kx=0$ with the underlying quadratic matrix polynomial
\begin{equation}\label{mm52}
P(\lambda):=\lambda^2M-\lambda(G-D)+K
\end{equation}
where $M,G,D,K\in\mathbb R^{n,n}$ satisfy $M=M^\top,D=D^\top,K=K^\top \geq0$ and $G=-G^\top$. As we will see below, this can be
viewed as a generalization of the dH structure to second order systems. The second  order case can be easily rewritten in first
order dH form but we will treat the problem directly in second order, and we will also discuss appropriate higher degree matrix
polynomials $P(\lambda) $ with an analogous structure.

Linear time invariant systems with the described  structure are very common in all areas of science and engineering
\cite{BeaMXZ18,MehM19,SchJ14} and typically arise via linearization arround a stationary solution. However, since all mathematical
models of physical systems are usually only
approximations of reality and data are typically inaccurate, it is an important question whether a given model is close to a model
with 'bad properties' such as an ill-posed  model without or with non-unique solution. To answer such questions for dH and
port-Hamiltonian systems has been an important research topic in recent years, see, e.g.,
\cite{AliMM20,BeaMV19,GilMS18,GilS17,GilS18,MehMS16,MehMS17}.

To classify whether a model is close to a 'bad model', one usually computes the distance to the nearest model with the 'bad' property.
In this paper we will discuss the distance to the set of singular matrix polynomials, i.e., those with a determinant that is identically zero,
and the distance to the nearest high-index problem, i.e., a problem with Jordan blocks associated to the eigenvalue $\infty$ of size
bigger than one.

While for general unstructured DAE systems the characterization of these two distances is very difficult and partially
open  \cite{BerGTWW17,BerGTWW19,ByeHM98,GugLM17,MehMW15}, the picture changes if one considers \emph{structured distances}, i.e.,
distances within the set of linear constant coefficient dHDAE systems. In this paper, we will make use of previous
results from \cite{GilMS18,MehMW18} to derive explicit characterizations for computing these distances in terms of null-spaces
of several matrices.

We use the following notation. By $\| X\|_F$ we denote the Frobenius norm
of a (possibly rectangular) matrix $X$, we extend this norm to matrix polynomials $P(\lambda)=\sum_{j=0}^k \lambda^j X_j$
by setting $\| P(\lambda)\|_F=\|[X_0,\dots,X_k]\|_F$. By $\lambda_{\min}(X)$ we denote the smallest eigenvalue of a positive
semidefinite matrix $X$.

The paper is organized as follows.
In Section~\ref{sec:prelim} we recall a few basic results about  linear time-invariant dH systems.
In Section~\ref{sec:prob} we present the different distances and state the main results for first order systems. Instead of immediately
presenting the corresponding proofs, we first consider related distance problems for a more general polynomial structure
in Section~\ref{sec:gendist}. These distance characterizations are then specialized in Section~\ref{sec:dhdaedist} to prove the main
results for the first order case. In Section~\ref{sec:qdh} we
consider corresponding distances for analogous quadratic matrix polynomials and also show how different representations of the first
order case can be related.

\section{Preliminaries}\label{sec:prelim}
We will make use of the Kronecker canonical form of a matrix pencil \cite{Gan59a}. Let us denote by
$\mathcal J_n(\lambda_0)$ the standard upper triangular Jordan block of size $n\times n$ associated with the eigenvalue $\lambda_0$
and let $\mathcal L_n$ denote the standard right Kronecker block of size $n\times(n+1)$, i.e.,
\[
\mathcal L_n=\lambda\left[\begin{array}{cccc}
1&0\\&\ddots&\ddots\\&&1&0
\end{array}\right]-\left[\begin{array}{cccc}
0&1\\&\ddots&\ddots\\&&0&1
\end{array}\right]\quad\mbox{and}\quad \mathcal J_n(\lambda_0)=\left[\begin{array}{cccc}
\lambda_0&1\\&\ddots&\ddots\\&&\ddots&1\\&&&\lambda_0\end{array}\right].
\]
\begin{theorem}[Kronecker canonical form]\label{th:kcf}
Let $E,A\in {\mathbb C}^{n,m}$. Then there exist nonsingular matrices
$S\in {\mathbb C}^{n,n}$ and $T\in {\mathbb C}^{m,m}$ such that
\begin{equation}\label{kcf}
S(\lambda E-A)T=\diag({\cal L}_{\epsilon_1},\ldots,{\cal L}_{\epsilon_p},
{\cal L}^\top_{\eta_1},\ldots,{\cal L}^\top_{\eta_q},
{\cal J}_{\rho_1}^{\lambda_1},\ldots,{\cal J}_{\rho_r}^{\lambda_r},{\cal N}_{\sigma_1},\ldots,
{\cal N}_{\sigma_s}),
\end{equation}
where $p,q,r,s,\epsilon_1,\dots,\epsilon_p,\eta_1,\dots,\eta_q,\rho_1,\dots,\rho_r,\sigma_1,\dots,\sigma_s\in\mathbb N$ and
$\lambda_1,\dots,\lambda_r\in\mathbb C$, as well as ${\cal J}_{\rho_i}^{\lambda_i}=I_{\rho_i}-\mathcal J_{\rho_i}(\lambda_i)$
for $i=1,\dots,r$ and $\mathcal N_{\sigma_j}=\mathcal J_{\sigma_j}(0)-I_{\sigma_j}$ for $j=1,\dots,s$.
This form is unique up to permutation of the blocks.
\end{theorem}
For real matrices (the case we discuss), a real version of the Kronecker canonical form is obtained under real transformation matrices $S,T$.
In this case the  blocks ${\cal J}_{\rho_j}^{\lambda_j}$ with $\lambda_j\in\Comp\setminus\Real$ have to be replaced with corresponding
blocks in real Jordan canonical form associated to the corresponding pair of conjugate complex eigenvalues, but the other
blocks have the same structure as in the complex case. An eigenvalue is called semisimple if the largest associated Jordan block has size one.

The sizes $\eta_j$ and $\epsilon_i$ of the rectangular blocks
are called the \emph{left and right minimal indices} of $\lambda E-A$, respectively.
The matrix pencil $\lambda E-A$, $E,A \in \mathbb C^{n,m}$ is called \emph{regular} if $n=m$ and
$\operatorname{det}(\lambda_0 E-A)\neq 0$ for some $\lambda_0 \in \mathbb C$,
otherwise it is called \emph{singular}.
A pencil is singular if and only if it has blocks of at least one of the types ${\cal L}_{\eps_j}$ or
${\cal L}^\top_{\eta_j}$ in the Kronecker canonical form.

The values $\lambda_1,\dots,\lambda_r\in\mathbb C$ are called the finite eigenvalues of $\lambda E-A$. If $s>0$, then
$\lambda_0=\infty$ is said to be an eigenvalue of $\lambda E-A$. (Equivalently, zero is then an eigenvalue of
the reversal $\lambda A-E$ of the pencil $\lambda E-A$.)
The sum of all sizes of blocks that are associated with a fixed eigenvalue
$\lambda_0\in\mathbb C\cup\{\infty\}$ is called the \emph{algebraic multiplicity} of $\lambda_0$.
The size of the largest block ${\cal N}_{\sigma_j}$ is
called the \emph{index} $\nu$ of the pencil $\lambda E-A$, where, by convention,  $\nu=0$ if $E$ is invertible.
The pencil is called {\em stable} if it is regular and if all eigenvalues are in the closed left half plane, and the ones lying on the
imaginary axis (including infinity) have the largest associated block of size at most one. Otherwise the pencil is called {\em unstable}.

The following result was shown in \cite{MehMW18}. We state the result in full generality, but clearly all statements also hold for the
special case that $E,Q,J,R$ are real and that $Q=I$ which is the case considered in this paper.
\begin{theorem}\label{thm:singind}
Let $E,Q\in \mathbb C^{n,m}$ satisfy $E^H Q=Q^H E\geq 0$ and let all left minimal indices of
$\lambda E-Q$ be equal to zero (if there are any). Furthermore, let
$J,R\in\mathbb R^{m,m}$ be such that we have $J=-J^H$, $R\geq 0$. Then the following
statements hold for the pencil $L(\lambda)=\lambda E-(J-R)Q$.
\begin{enumerate}
\item[\rm (i)] If $\lambda_0\in\mathbb C$ is an eigenvalue of $L(\lambda)$ then $\operatorname{Re}(\lambda_0)\leq 0$.
\item[\rm (ii)] If $\omega\in\mathbb R\setminus\{0\}$ and $\lambda_0=i\omega$ is an eigenvalue of $L(\lambda)$, then
$\lambda_0$ is semisimple. Moreover, if the columns of $V\in\mathbb C^{m,k}$ form a basis of a regular deflating
subspace of $L(\lambda)$ associated with $\lambda_0$, then $RQV=0$.

If, additionally, $Q$ is nonsingular then the
previous statement holds for $\lambda_0=0$ as well. If $Q$ is singular then $\lambda_0=0$ need not be semisimple,
but if $L(\lambda)$ is regular, then Jordan blocks associated with $\lambda_0=0$ have size at most two.
\item[\rm (iii)] The index of $L(\lambda)$ is at most two.
\item[\rm (iv)] All right minimal indices of $L(\lambda)$ are at most one (if there are any).
\item[\rm (v)] If in addition $\lambda E-Q$ is regular, then all left minimal indices of $L(\lambda)$ are zero (if there are any).
\end{enumerate}
\end{theorem}
\proof
For the proof see \cite{MehMW18}. The additional statement in (ii) on the eigenvalue $\lambda_0=0$ was not presented in \cite{MehMW18},
but it follows in a straightforward manner from \cite[Theorem~6.1]{MehMW18} and the proof of \cite[Corollary~6.2]{MehMW18}.
\eproof

Theorem~\ref{thm:singind} illustrates that the special structure of dH systems imposes many restrictions in the
spectral data and this has also an advantage when determining the distances to the nearest 'bad' problem. In particular,
Theorem~\ref{thm:singind} implies that the distance to instability and the distance to higher index coincide for a pencil $L(\lambda)$
with $Q$ nonsingular.

The following well-known lemma, see \cite{BreCP96,KunM06} (also stated for the general complex case), will be needed in order
to make statements about the index of a matrix pencil in special situations.
\begin{lemma}\label{lem:index}
Let $E,A\in\mathbb C^{n,n}$ be matrices of the form
\[
E=\mat{cc}E_{11}&0\\ 0&0\rix\quad\mbox{and}\quad A=\mat{cc}A_{11}&A_{12}\\ A_{21}&A_{22}\rix,
\]
where $E_{11}$ is invertible.
\begin{enumerate}[\rm (i)]
\item If $A_{22}$ is invertible, then the pencil $\lambda E-A$ is regular and has index one;
\item if $A_{22}$ is singular, then the pencil $\lambda E-A$ is singular or has an index greater than or equal to two.
\end{enumerate}
\end{lemma}

\section{Problem statement and main results for dHDAE systems}\label{sec:prob}

We are interested in the following distance problems for matrix pencils $L(\lambda)$ of the form \eqref{dHwoQ}
under perturbations that preserve the special structure of the pencil.
%
\begin{definition}\label{def:distances}
Let $\mathcal{L}$ denote the class of square $n\times n$ real matrix pencils of the form \eqref{dHwoQ}.
%
%
Then
\begin{enumerate}[1)]
\item  the  \emph{structured distance to singularity}  is defined as
\begin{equation}\label{distEJR}
d_{\sing}^{\mathcal L}\big(L(\lambda)):=\inf\big\{\big\|  \Delta_L(\lambda) \big\|_F\ \big|\
L(\lambda)+\Delta_L(\lambda)\in\mathcal{L} \mbox{ and is singular}\big\};
\end{equation}
\item the  \emph{structured distance to the nearest high-index problem} is defined as
\begin{equation}\label{indexdistEJR}
d_{\hi}^{\mathcal{L}}\big(L(\lambda)):=\inf\big\{
\big\|  \Delta_L(\lambda)
\big\|_F\ \big| \
L(\lambda)+\Delta_L(\lambda)\in \mathcal{L} \mbox{ and is of index}\geq 2\big\};
\end{equation}
\item the  \emph{structured distance to instability} is defined as
\begin{equation}\label{instdistEJR}
d_{\inst}^{\mathcal{L}}\big(L(\lambda)\big):=\inf\big\{
\big\|
\Delta_L(\lambda)
\big\|_F \
\big| \ L(\lambda)+\Delta_L(\lambda)\in\mathcal{L} \mbox{ and is unstable}
\big \}.
\end{equation}
%
\end{enumerate}
\end{definition}
Note that all defined distances are meaningful, as for each matrix $X\in\Real^{n,n}$ the decomposition into a sum $X=X_1+X_2$ of a
skew-symmetric matrix
$X_1=\frac12(X-X^\top)$  and  symmetric  matrix $X_2=\frac12(X+X^\top)$  is unique. Furthermore, we have
$\norm X_F^2= \norm{ X_1}_F^2 + \norm{ X_2}_F^2=\big\|[X_1,X_2]\big\|_F^2$
due to the trace of $X_1^\top X_2$ being zero.
Thus, the constraint $L(\lambda)+\Delta_L(\lambda)\in\mathcal{L}$ in \eqref{distEJR}--\eqref{instdistEJR} is the same as writing
\[
\Delta_L(\lambda)=\lambda \Delta_E -(\Delta_J - \Delta_R),
\]
with $ \Delta_J=-\Delta_J^\top$ and $E+\Delta_E,R+\Delta_R\geq 0$, and we have
$\big\| [ \Delta_J, \Delta_R, \Delta_E]\big\|_F=\big\| [ \Delta_L(\lambda)] \big\|_F$.
The positivity conditions for $E+\Delta_E,R+\Delta_R$ are crucial.
Examples presented in Section \ref{sec:strunstr} show that they can neither be omitted nor simplified to
$E+\Delta_E,R+\Delta_R$ being merely symmetric.
%
\begin{theorem}\label{Q=I}
Let  $L(\lambda)=\lambda E -(J-R)\in \mathcal L$.
Then the following statements hold.
\begin{enumerate}[\rm (i)]
\item\label{I4} The pencil $L(\lambda)$ is singular if and only if $\,\ker J\cap\ker E\cap\ker R\neq\set0$. In that case there exists
an orthogonal transformation matrix $U\in\Real^{n,n}$ such that
\[
U^\top EU=\mat{cc} E_{11} & 0 \\ 0 & 0 \rix, \quad  U^\top JU=\mat{cc} J_{11} & 0 \\ 0 & 0	 \rix,\quad U^\top RU=
\mat{cc} R_{11} & 0 \\ 0 & 0 \rix,
\]
where the pencil $\lambda E_{11} - (J_{11}-R_{11})$ is regular and has the size $(n-r)\times(n-r)$ with $r=\dim\ker(E-J+R)>0$.
In particular, all right and left minimal indices of $L(\lambda)$ in its Kronecker canonical form are zero.
\item\label{I3} The index of $L(\lambda)$ is at most two. Furthermore, the following statements are equivalent.
\begin{enumerate}[\rm (a)]
\item For any $\eps>0$ there exists a
pencil $\lambda \widetilde E- (\widetilde J-\widetilde R) $ with $\widetilde E,\widetilde R\geq 0$, and $\widetilde J=-\widetilde J^\top$
which is regular and of index two such that

\begin{equation}\label{tildesmall}
\big\|\big[E-\widetilde E,\ J-\widetilde J,\ R-\widetilde R\big]\big\|_F
=\big\|\big[E-\widetilde E, (J-R)-(\widetilde J-\widetilde R)\big]\big\|_F\leq \eps,
\end{equation}

i.e., $L(\lambda)$ is in the closure of the set of regular dH pencils of index two.
\item $\ker E\cap\ker R\neq\set0$.
\end{enumerate}
\end{enumerate}
\end{theorem}
%

To construct the perturbations where the distance to singularity ia achieved, we use the following ansatz. For a matrix $Y\in\mathbb R^{n,n}$
and a vector $u\in\mathbb R^n$ with $\|u\|_2=1$ we  define the matrix
\begin{equation}\label{deltaY}
\Delta_Y^u=-uu^\top Y-Yuu^\top+uu^\top Yuu^\top,
\end{equation}
that will be used at several occasions during the paper.
Then we obtain the  following characterization of the distance to singularity.
\begin{theorem}\label{singind}
Let  $\lambda E -(J-R) \in \mathcal L$.
Then the following statements hold.
\begin{enumerate}[\rm (i)]
\item\label{singindI}
The distance to singularity \eqref{distEJR} is attained with a perturbation
$\Delta_E=\Delta_E^u$, $\Delta_J=\Delta_J^u$, and $\Delta_R=\Delta_R^u$ as in~\eqref{deltaY}
for some $u\in\mathbb R^n$ with $\|u\|_2=1$. The distance is given by
%
\begin{multline*}
d_{\sing}^{\mathcal{L}}\big(\lambda E - (J-R)\big)\\
=\min_{u \in\Real^n\atop\norm u=1} \sqrt{2\norm{Ju}^2+2\big\|(I-uu^\top)Eu\big\|^2+(u^\top Eu)^2+
2\big\|(I-uu^\top)Ru\big\|^2+(u^\top Ru)^2}
\end{multline*}
%
and is bounded as
\begin{equation}\label{lmino}
\sqrt{\lambda_{\min} (-J^2+R^2+E^2)}\leq d_{\sing}^{\mathcal{L}}\big(\lambda E - (J-R)\big)\leq  \sqrt{2\cdot\lambda_{\min}(-J^2+R^2+E^2)}.
\end{equation}
\item\label{singindII} The distance to higher index~\eqref{indexdistEJR} and the distance to instability~\eqref{instdistEJR} coincide and
satisfy
\begin{multline*}
d_{\hi}^{\mathcal{L}}\big(\lambda E - (J-R)\big)=d_{\inst}^{\mathcal{L}}\big(\lambda E - (J-R)\big)\\
=\min_{u \in\Real^n\atop\norm u=1} \sqrt{2\big\|(I-uu^\top)Eu\big\|^2+(u^\top Eu)^2+
2\big\|(I-uu^\top)Ru\big\|^2+(u^\top Ru)^2} 
\end{multline*}
and are bounded as
\begin{equation*}
\sqrt{\lambda_{\min}(E^2+R^2)}\leq d_{\hi}^{\mathcal{L}}\big(\lambda E - (J-R)\big)=d_{\inst}^{\mathcal{L}}\big(\lambda E - (J-R)\big)
\leq  \sqrt{2\cdot\lambda_{\min}(E^2+R^2)}.
\end{equation*}
\end{enumerate}
\end{theorem}
The proofs of Theorems~\ref{Q=I} and~\ref{singind} are given in Section~\ref{s:proof}, where they are obtained as
simple consequences of a general theory developed in Section \ref{s:distpoly} for matrix polynomials with a special symmetry structure.
Before  we give the proofs, we will first consider a more general minimization problem in the next section.
\section{General distance problems}\label{sec:gendist}
In this section, we present a solution to a quite general minimization problem. This will  allow us to solve
the distance problems for dH pencils introduced in Section~\ref{sec:prob} as well as analogous problems for structured matrix polynomials
with a dH like structure in a unified manner.

Theorem~\ref{Q=I} states that both the distance to singularity as well as to higher index for a dH pencil
as in~\eqref{dHwoQ} can be expressed via the existence of a common kernel of two or three structured matrices, so that both problems
can be reinterpreted as a distance problem to the common kernel of matrices with symmetry and positivity structures.
This concept will now be extended to more than three matrices.

\subsection{Distance to the common kernel of a tuple of structured matrices}\label{subsec:ck}

\begin{definition}\rm Let $\mathcal S^n_\ell$ denote the following set of $(\ell+2)$-tuples  of $n\times n$ real matrices
\[
\mathcal{S}^n_\ell:=\big\{(J,X_0,\dots,X_\ell)\in(\mathbb R^{n,n})^{\ell+2}\,\big|\,
J^\top=-J,\ {X_i}=X_i^T\geq 0, \ i=1,\dots,\ell\big\},
\]
where $\ell\geq 0$ and $n\geq 1$ are fixed.
For a given tuple $(J,X_0,\dots,X_{\ell})\in\mathcal S^n_\ell$ we define the {\em structured distance to the common kernel}
$d_{\ck}^{\mathcal S^n_\ell}(J,X_0,\dots,X_\ell )$
as
\begin{equation}\label{distck}
\inf\left\{\big\| \left[\Delta_J,\Delta_{X_0},\dots,\Delta_{X_\ell }\right]\big\|_F\Big|\
\begin{array}{c} (J+\Delta_J,X_0+\Delta_{X_0},\dots,X_\ell+\Delta_{X_\ell})\in\mathcal S^n_\ell,\\ \ker(J+\Delta_J)
\cap\bigcap_{i=0}^\ell \ker(X_i+\Delta_{X_i})\neq\{0\}
\end{array} \right\}.
\end{equation}
In the following, we often drop the dependence on $\ell$ and $n$ in the notation for simplicity, thus writing
$d_{\ck}^{\mathcal S}(J,X_0,\dots,X_\ell)$.
\end{definition}
Observe that in determining $d_{\ck}^{\mathcal S}(J,X_0,\dots,X_\ell )$ we  measure the distance to a closed set.
\begin{lemma}\label{closedJ}
The set of all $(J,X_0,\dots,X_\ell)\in\mathcal S^n_\ell$ satisfying
$\ker J\cap \ker X_0\cap \dots\cap\ker X_\ell \neq \set 0$ is a closed subset in $\mathcal S^n_\ell$.
\end{lemma}
\proof
The proof follows by  considering sequences of tuples $(J^{(m)},X_0^{(m)},\dots,X_\ell^{(m)})$ and a convergent subsequence of
a sequence of unit vectors $u_m$ satisfying
\[
J^{(m)}u_m = X_i^{(m)}u_m =0,\quad i=1,\dots\ell.\quad\mbox{\eproof}
\]
Before we present the solution of the minimization problem, we first develop equivalent conditions for $J,X_0,\dots,X_\ell$
to have a nontrivial common kernel.
\begin{proposition}\label{charsinggen-new}
Let $(J,X_0,\dots,X_\ell) \in\mathcal S^n_\ell$. Then
%
\begin{equation}\label{6.10.19a}
\ker J\cap \ker X_0\cap \dots\ker X_\ell = \ker (J^\top J+X_0^2+\cdots+X_\ell ^2)= \ker (-J+X_0+\cdots+X_\ell).
\end{equation}
Furthermore, there exists an orthogonal matrix $U\in\mathbb R^{n,n}$ such that
\begin{equation}\label{e1j1r1}
U^\top JU=\mat{cc}\widetilde J&0\\ 0&0\rix,\quad U^\top X_iU=\mat{cc}\widetilde X_i&0\\ 0&0\rix,\quad i=0,\dots,\ell
\end{equation}
with some $\widetilde J,\widetilde X_0,\dots,\widetilde X_\ell \in\mathbb R^{n-r,n-r}$, where $r=\dim \ker (-J+X_0+\cdots+X_\ell)\geq 0$ and
where the matrix $-\widetilde J+\widetilde X_0+\cdots+\widetilde X_\ell$ is invertible.
\end{proposition}
\proof
The inclusion $\ker J\cap \ker X_0\cap \dots\ker X_\ell \subseteq \ker (J^\top J+X_0^2+\cdots+X_\ell ^2)$ is trivial.
To prove the converse, let $x\in\ker(J^TJ+X_0^2+\dots+X_\ell ^2)$ be nonzero. Since each summand is positive semidefinite, we
obtain $X_0^2x=\cdots=X_\ell ^2x=0$ and $J^2x=-J^TJx=0$.
Noting that $\ker Y=\ker Y^2$ holds for any symmetric or skew-symmetric matrix finishes the proof.

The inclusion $\ker J\cap \ker X_0\cap \dots\ker X_\ell \subseteq \ker (-J+X_0+\cdots+X_\ell)$ is again trivial. To prove the converse
let $x\in\ker(-J+X_0+\dots+X_\ell )$ be nonzero. Since $x^\top Jx=0$, we obtain that $x^\top X_0x+\cdots+x^\top X_\ell x=0$ and
since each of the matrices $X_0,\dots,X_\ell$ is positive semidefinite, we obtain $X_0x=\cdots=X_\ell x=0$, which then implies $Jx=0$ as well.

To prove the last assertion, let $U$ be an orthogonal matrix with last $r$ columns spanning the kernel  of $-J+X_0+\cdots+X_\ell$.  Note that
if $u$ is one of those $r$ last columns of $U$ then~\eqref{6.10.19a} implies that $u^\top J=0$ and $u^\top X_i=0$ for $i=0,\dots,\ell$,
which shows the formula \eqref{e1j1r1}.
\eproof
\begin{remark}\label{rem:nono}{\rm
We highlight that the nonegativity assumption for the matrices $X_i$ is crucial for the two nontrivial inclusions in
Proposition~\ref{charsinggen-new}. For example, consider
\[
J=\mat{cc}0&1\\ -1&0\rix\quad\mbox{and}\quad X_0=\mat{cc}1&0\\ 0&-1\rix.
\]
Then $J-X_0$ is singular while the intersection of the kernels of $J$ and $X_0$ is trivial.

Also note that while an arbitrarily large number of symmetric positive semidefinite matrices $X_0,\dots,X_\ell$
can be considered, the results from Proposition~\ref{charsinggen-new} are no longer true if a second skew-symmetric matrix is
involved. For example, consider the matrices
\[
J_1=\mat{ccc}0&1&0\\ -1&0&0\\ 0&0&0\rix\quad\mbox{and}\quad J_2=\mat{ccc}0&0&1\\ 0&0&0\\ -1&0&0\rix
\]
Then $J_1+J_2$ is singular (in fact, even the pencil $\lambda J_1+J_2$ is singular), but $J_1$ and $J_2$ do not have a common kernel.}
\end{remark}

Given $(J,X_0,\dots,X_\ell) \in\mathcal S^n_\ell$ as in Proposition~\ref{charsinggen-new}, we aim to characterize all perturbations
that produce a nontrivial
common kernel of the matrices $J,X_0,\dots,X_\ell$ while preserving their individual structures. For this, we will use particular
perturbations whose special properties will be presented in the following lemma.
%
%
%
\begin{lemma}\label{lem:10.9.19}
Let $Y\in \Real^{n,n}$, let $u\in \Real^{n,n}$ be a vector with $\|u\|_2=1$ and let
\begin{equation}\label{derpert}
\Delta^u_Y:=-uu^\top Y-Yuu^\top+uu^\top Y uu^\top.
\end{equation}
Then the following statements hold.
\begin{enumerate}[\rm (i)]
\item\label{l1}$u\in\ker(Y+\Delta_Y^u)$, in particular, $Y+\Delta_Y^u$ is singular.
\item\label{l2} $\rank\ \Delta_Y^u\leq 2$, and $\rank\ \Delta_Y^u\leq 1$ if and only if $u$ is a
right or left eigenvector of $Y$.
\item\label{l3} $\|\Delta^u_Y\|^2_F=\big\|(I-uu^\top)Yu\big\|^2_F+\big\|u^\top Y(I-uu^\top)\big\|_F^2+(u^\top Yu)^2$.
\item\label{l4} If $Y\geq 0$ then $Y+\Delta_Y^u\geq 0$ and $\|\Delta^u_Y\|^2_F=2\big\|(I-uu^\top)Yu\big\|^2_2+(u^\top Yu)^2$.
\item\label{l5} If $Y^\top=-Y$, then $\Delta^u_Y=-uu^\top Y-Yuu^\top$ and $\|\Delta^u_Y\|^2_F=2\|Yu\|_2$.
\end{enumerate}
\end{lemma}
\proof
\eqref{l1} immediately follows from $\Delta^u_Yu=-Yu$. For the proof
of \eqref{l2} let $U\in\mathbb R^{n,n}$ be an orthogonal matrix with last column $u$. Then we obtain
\begin{equation}\label{eq:10.9.19}
U^TYU=\mat{cc}Y_{11}&Y_{12}\\ Y_{21}&Y_{22}\rix\quad\mbox{and}\quad U^T\Delta_Y^uU=\mat{cc}0&-Y_{12}\\ -Y_{21}&-Y_{22}\rix
\end{equation}
for some $Y_{11},Y_{12},Y_{21},Y_{22}$ with $Y_{11}\in\mathbb R^{n-1,n-1}$ which immediately shows that $\rank\ \Delta_Y^u\leq 2$.
In particular, we have $\rank \Delta_Y^u\leq 1$ if and only if $Y_{12}=0$ or $Y_{21}=0$ which is equivalent to $u$ being a right or left
eigenvector of $Y$, respectively. Moreover, \eqref{l3} immediately follows from the representation~\eqref{eq:10.9.19} using
that
\[
(I-uu^\top)Yu=\mat{c}Y_{12}\\ 0\rix,\quad uY(I-uu^\top)=\mat{cc}Y_{21}&0\rix,\quad\mbox{and}\quad Y_{22}=u^\top Yu.
\]
Finally, using the additional (skew-)symmetry structure, we obtain \eqref{l4} and \eqref{l5}, where the part
$Y+\Delta_Y^u\geq 0$ in \eqref{l4} again follows from the representation~\eqref{eq:10.9.19}.
\eproof

We highlight that the first property of statement~\eqref{l4} in Lemma~\ref{lem:10.9.19} will become essential in what follows, because
it allows us to perform
a perturbation that makes a symmetric matrix singular while simultaneously preserving the positive semidefiniteness of the matrix.
With these preparations, we obtain the following theorem that characterizes structure-preserving
perturbations to matrices with a nontrivial common kernel.
\begin{theorem}\label{thm:smallper}
Let $(J,X_0,\dots,X_\ell) \in\mathcal S^n_\ell$, i.e., $J^\top=-J$ and $X_i^\top=X_i\geq 0$ for $i=0,\dots,\ell$. Furthermore,
for any $u\in\Real^n$, $\|u\|_2=1$, consider the perturbation matrices
\begin{equation}\label{derpertold}
\Delta^u_J:= -uu^\top J - J uu^\top\quad\mbox{and}\quad
\Delta^u_{X_i}:= -uu^\top X_i - X_i uu^\top+ uu^\top X_i uu^\top ,\quad i=0,\dots,\ell .
\end{equation}
Then the following statements hold.
\begin{enumerate}[\rm (i)]
\item\label{t1} For any vector $u\in\Real^n$, $\|u\|_2=1$, we have
\begin{equation}\label{deltaprop}
(\Delta^u_J)^\top=-\Delta^u_J,\quad\mbox{as well as}\quad(\Delta^u_{X_i})^\top=\Delta_{X_i}^u\;\mbox{and}\;
X_i+\Delta_{X_i}^u\geq 0, \quad i=0,\dots,\ell .
\end{equation}
Furthermore, the kernels of the matrices $J+\Delta^u_J$, $X_0+\Delta_{X_0}^u,\dots,X_\ell +\Delta_{X_\ell }^u$ have a nontrivial intersection.
\item\label{t4} For any vector $u\in\Real^n$, $\|u\|_2=1$, we have
\[
\norm{\Delta^u_J}^2_F=2\norm{Ju}_2,\quad\mbox{and}\quad \norm{\Delta^u_{X_i}}^2_F= 2\big\|(I-uu^\top) X_iu\big\|^2+(u^\top X_iu)^2,
\quad i=0,\dots,\ell .
\]
\item\label{t3}  Let $\Delta_J,\Delta_{X_0},\dots,\Delta_{X_\ell }\in\mathbb R^{n,n}$ be any perturbation matrices satisfying
\begin{equation}\label{deltaprop2}
\Delta_J^\top=-\Delta_J\quad\mbox{as well as}\quad \Delta_{X_i}^\top=\Delta_{X_i},\;\mbox{and}\;X_i+\Delta_{X_i}\geq 0,\quad i=0,\dots,\ell,
\end{equation}
and such that the kernels of the matrices $J+\Delta_J$, $X_0+\Delta_{X_0},\dots,X_\ell +\Delta_{X_\ell }$ have a nontrivial intersection. Then
\[
\norm{\Delta_J^u}_F\leq \norm{ \Delta_J}_F,\quad\mbox{and}\quad \norm{\Delta_{X_i}^u}_F\leq \norm{ \Delta_{X_i}}_F,\quad i=0,\dots,\ell .
\]
for some real vector $u$ with $\|u\|_2=1$
\end{enumerate}
\end{theorem}
\proof
\eqref{t1} and~\eqref{t4} follow immediately from Lemma~\ref{lem:10.9.19}.
To prove \eqref{t3}, consider any perturbation matrices $\Delta_J,\Delta_{X_0},\dots,\Delta_{X_\ell }$	
satisfying \eqref{deltaprop2} such that the kernels of the matrices $J+\Delta_J$, $X_0+\Delta_{X_0},\dots,X_\ell +\Delta_{X_\ell }$
have a nontrivial intersection. Then by Proposition \ref{charsinggen-new}, there exists an orthogonal matrix $U$ such that
\begin{equation}\label{eq:10.9.19a}
U^\top (J+\Delta_J)U=\mat{cc}\widetilde J&0\\ 0&0\rix,\quad U^\top (X_i+\Delta_{X_i})U=\mat{cc}\widetilde X_i&0\\ 0&0\rix,\quad i=0,\dots,\ell
\end{equation}
with some $\widetilde J,\widetilde X_0,\dots,\widetilde X_\ell \in\mathbb R^{n-1,n-1}$, not necessarily invertible, i.e., in contrast
to~\eqref{e1j1r1} we split only one vector from the intersection of kernels. Transforming and decomposing accordingly, we have
\begin{equation}\label{eq:10.9.19b}
U^\top J U=\matp{ \widetilde K  & t \\ -t^\top & 0},\quad\mbox{and}\quad
U^\top X_i U=\matp{ \widetilde S_i  & s_i \\ s_i^\top  & r_{i}},\quad i=0,\dots,\ell
\end{equation}
for some skew-symmetric matrix $\widetilde K\in\mathbb R^{n-1,n-1}$, some symmetric matrices $\widetilde S_i\in\mathbb R^{n-1,n-1}$,
some $r_{i}\in\Real$, $s_i\in\Real^{n-1}$  ($i=0,\dots,\ell$), and some $t\in\Real^{n-1}$.
Subtracting~\eqref{eq:10.9.19b} from~\eqref{eq:10.9.19a}, we obtain that
\begin{equation}
U^\top\Delta_JU=\matp{ \widetilde J-\widetilde K & -t  \\ t^\top & 0 },\quad\mbox{and}\quad
U^\top\Delta_{X_i}U=\matp{ \widetilde X_i-\widetilde S_i& -s_i \\ -s_i^\top  & -r_{i} },\quad i=0,\dots,\ell .
\end{equation}
Observe that for the particular choice $u=Ue_n$ the perturbations from~\eqref{derpert} have, by \eqref{eq:10.9.19b}, the forms
\begin{equation}
U^\top\Delta_J^uU=-\matp{ 0 & -t  \\ t^\top & 0}\quad\mbox{and}\quad
U^\top\Delta_{X_i}^uU=-\matp{0 & s_i \\ s_i^\top  & r_{i} },\quad i=0,\dots,\ell .
\end{equation}
%
Since the Frobenius norm is invariant under real orthogonal transformations, we immediately obtain
that $\norm{\Delta_J^u}_F\leq \norm{ \Delta_J}_F$ and $\norm{\Delta_{X_i}^u}_F\leq \norm{ \Delta_{X_i}}_F$ for $i=0,\dots,\ell$.
\eproof
We now have all ingredients to state and prove the solution of our general minimization problem.
\begin{theorem}\label{main}
Let $(J,X_0,\dots,X_\ell) \in\mathcal S^n_\ell$, i.e., $J^\top=-J$ and $X_j^\top=X_j\geq 0$ for $j=1,\dots,\ell$.
Then the structured distance $d_{\ck}^{\mathcal S}(J,X_0,\dots,X_\ell )$ to the common kernel~\eqref{distck} is attained at
 $\Delta_J=\Delta_J^u$, $\Delta_{X_0}=\Delta_{X_0}^u,\dots,\Delta_{X_\ell }=\Delta_{X_\ell }^u$ being as in \eqref{derpertold}
for some $u\in\mathbb R^{n}$ with $\|u\|_2=1$. Consequently,
\begin{equation*}
d_{\ck}^{\mathcal S}(J,X_0,\dots,X_\ell )=\min_{ u \in\Real^n, \norm u=1} \left(2\norm{Ju}_2+\sum_{i=1}^\ell \Big(
2\big\|(I-uu^\top) X_iu\big\|^2+(u^\top X_iu)^2\Big)\right)^{1/2},
\end{equation*}
and in addition, we have the bounds
\begin{equation}\label{lmin}
\sqrt{\lambda_{\min} (-J^2+X_0^2+\cdots+X_\ell ^2)}\leq d_{\ck}^{\mathcal S}(J,X_0,\dots,X_\ell )\leq \sqrt{2\cdot \lambda_{\min}
(-J^2+X_0^2+\cdots+X_\ell ^2)}.
\end{equation}
\end{theorem}
\proof
The first two statements follow directly from Theorem \ref{thm:smallper}. It remains to prove \eqref{lmin}. For this aim note that
for every $u\in \Real^n$ with $\|u\|_2=1$, we have
\[
\begin{split}
u^\top (-J^2+X_0^2+\cdots+X_\ell ^2) u  =&  \norm{Ju}^2 +\norm{X_0u}^2+\cdots+\norm{X_\ell u}^2 \\
=& \left(\norm{Ju}^2+\sum_{i=1}^\ell \Big(\norm{(I-uu^\top)X_iu}^2 +(u^\top X_iu)^2\Big)  \right).
\end{split}
\]
Taking the infimum over all $u\in \Real^n$ with $\|u\|_2=1$ shows \eqref{lmin}.
\eproof
\begin{remark}{\rm
In the special case $\lambda_{\min} (-J^2+X_0^2+\cdots+X_\ell ^2)=0$ it immediately follows that
$d_{\ck}^{\mathcal{S}}(J,X_0,\dots,X_\ell )=0$. This is in line with Proposition~\ref{charsinggen-new}, because the singularity of the
matrix $-J^2+X_0^2+\cdots+X_\ell ^2$ is equivalent to the existence of a nontrivial common kernel of the matrices
$J,X_0,\dots,X_\ell$.
}
\end{remark}

\subsection{Distance problems for structured matrix polynomials}\label{s:distpoly}

As a first application of the results from Subsection~\ref{subsec:ck}, we will consider distance problems for
a particular class of structured matrix polynomials.
To this end, recall from \cite{Gan59a} that by definition a square matrix polynomial $P(\lambda)=\sum_{i=0}^k \lambda^i Y_i$ is
\emph{singular} if and only if $\det P(\lambda)\equiv 0$. Also recall that the
\emph{companion linearization}
\begin{equation}\label{companion}
L(\lambda)=\lambda\mat{cccc} Y_k&&&\\ &I&&\\ &&\ddots &\\ &&&I\rix+\mat{cccc}Y_{k-1}&\dots&Y_1&Y_0\\ -I&0&&\\ &\ddots&\ddots&\\ &&-I&0\rix.
\end{equation}
of $P(\lambda)$ is a \emph{strong linearisation} in the sense of \cite{DeDM09a}. In particular, $L(\lambda)$ is singular if and only
if $P(\lambda)$ is singular. Furthermore, as shown in \cite{DeDM09a}, in the linearization the spectral data
for eigenvalues of $P(\lambda)$ is preserved. Therefore,
  for the sake of simplicity, we define the notions of {\em index} and {\em instability} for the  matrix polynomial $P(\lambda)$ via
  the respective notions of the Kronecker canonical form of~\eqref{companion}, cf. Section \ref{sec:prelim}. We then extend
  Definition \ref{def:distances} as follows.
\begin{definition}
Consider the class of  matrix polynomials
\[
\mathcal P^n_{k,j}:=\left.\left\{-\lambda^j J+\sum_{i=0}^ k \lambda^iA_i\,\right|\,J^\top=-J,\
A_i^\top=A_i\geq 0 \in\mathbb R^{n,n}, \ i=0,\dots,k\right\},
\]
where $n\geq 1$, $k,j \geq 0$ and, without loss of generality, $j \leq k$. Then for $P(\lambda)\in\mathcal P^n_{k,j}$
\begin{enumerate}[1)]
\item  the  \emph{structured distance to singularity}  is defined as
\begin{equation}\label{distP}
d_{\sing}^{\mathcal P^n_{k,j}}\big(P(\lambda)):=\inf\big\{\big\|  \Delta_P(\lambda) \big\|_F\ \big|\
P(\lambda)+\Delta_P(\lambda)\in\mathcal P^n_{k,j} \mbox{ is singular}\big\};
\end{equation}
\item the  \emph{structured distance to the nearest high index problem} is defined as
\begin{equation}\label{indexdistP}
d_{\hi}^{\mathcal P^n_{k,j}}\big(P(\lambda)):=\inf\big\{\big\|  \Delta_P(\lambda) \big\|_F\ \big|\
P(\lambda)+\Delta_P(\lambda)\in\mathcal P^n_{k,j} \mbox{ is of index}\geq 2\big\};
\end{equation}
\item the  \emph{structured distance to instability} is defined as
\begin{equation}\label{instdistP}
d_{\inst}^{\mathcal P^n_{k,j}}\big(P(\lambda)\big):=\inf\big\{\big\|  \Delta_P(\lambda) \big\|_F\ \big|\
P(\lambda)+\Delta_P(\lambda)\in\mathcal P^n_{k,j}\mbox{ is unstable}\big\}.
\end{equation}
\end{enumerate}
We often simply write $d_{\star}^{\mathcal P}\big(P(\lambda))$ instead of
$d_{\star}^{\mathcal P^n_{k,j}}\big(P(\lambda))$ for $\star\in\{\sing,\hi,\inst\}$.
\end{definition}
In other words, $\mathcal P^n_{k,j}$ consists of the set of matrix polynomials of degree less than or equal to $k$
for which all coefficients are symmetric positive semidefinite except for the coefficient at $\lambda^j$ which
is only assumed to have a positive semidefinite symmetric part. Particular examples for this kind of matrix polynomials
are the dH pencils of the form~\eqref{dHwoQ}, i.e., the set $\mathcal P^n_{1,0}$, and quadratic matrix polynomials of
the form~\eqref{mm52}, i.e., the set $\mathcal P^n_{2,1}$.
 Observe that if both $P(\lambda),P(\lambda)+\Delta_P(\lambda)\in\mathcal P^n_{k,j}$ then $\Delta_P(\lambda)$ must take the form
\[
\Delta_P(\lambda)=-\lambda^j \Delta_J+\sum_{i=0}^k\lambda^i\Delta_{A_i},
\]
where $\Delta_J^\top=-\Delta_J$ and $\Delta_{A_i}^\top=\Delta_{A_i}$.
We have the following theorem for characterizing the distance to the nearest singular or high index matrix polynomial.
\begin{theorem}\label{distances}
Let $k\geq 1$ and $j\in\{0,\dots,k\}$ and consider the set $P^n_{k,j}$ of matrix polynomials
\[
P(\lambda)=-\lambda^j J+\sum_{i=0}^k\lambda^iA_i.
\]
with $J,A_0,\dots,A_k\in\mathbb R^{n,n}$, $J^\top=-J$ and $A_i^\top=A_i\geq 0$ for $i=0,\dots,k$.
\begin{enumerate}[\rm (i)]
\item\label{tg1} If $P(\lambda)\in\mathcal P^n_{k,j}$ then the following statements are equivalent:
\begin{enumerate}[\rm (a)]
\item\label{p1} the polynomial $P(\lambda)$ is singular, i.e., $\det P(\lambda)\equiv 0$;
\item\label{p2} the matrix $P(1)$ is singular;
\item\label{p3} the kernels of the matrices $J,A_0,\dots,A_k$ have a nontrivial intersection.
\end{enumerate}
\item\label{tg1a} If $P(\lambda)\in\mathcal P^n_{k,j}$, then  its distance to the set of singular matrix
polynomials in $\mathcal P^n_{k,j}$ equals the distance to the common kernel of the matrices $J,A_0,\dots,A_k$, i.e.,
\begin{equation}\label{distdist1}
d_{\sing}^{\mathcal P}(P(\lambda))=d_{\ck}^{\mathcal S}(J,A_0,\dots,A_k).
\end{equation}
\item\label{tg2} If $\max\{n,k\}>1$, then the closure of the set
\begin{equation}\label{eins}
\mathcal I_{\hi}:=\big\{P(\lambda)\in\mathcal P^n_{k,j}\,\big|\, P(\lambda)\;\mbox{\rm is regular and has index greater than one}\big\}
\end{equation}
in $\mathcal P^n_{k,j}$ is equal to
\begin{equation}\label{zwei}
\mathcal K:=\left.\left\{-\lambda^j J+\sum_{i=0}^k\lambda^iA_i\;\right|\;\ker A_k\cap\ker A_{k-1}\neq\{0\} \right\}\quad \mbox{if } j<k
\end{equation}
or to
\[
\mathcal K:=\left.\left\{-\lambda^k J+\sum_{i=0}^k\lambda^iA_i\;\right|\;\ker J\cap\ker A_k\cap\ker A_{k-1}\neq\{0\} \right\}
\quad\mbox{if } j=k>1.
\]
If $\max\{n,k\}=1$ or $j=k=1$, then $\mathcal I_{\hi}$ is empty.
\item\label{tg2a} Let $P(\lambda)\in\mathcal P^n_{k,j}$. If $\max\{n,k\}>1$, then the distance of $P(\lambda)$ to the set of higher index
polynomials in
$\mathcal P^n_{k,j}$ equals the distance to the respective common kernel
\begin{equation}\label{distdist2}
d_{\hi}^{\mathcal P}(P(\lambda))=\begin{cases}  d^{\mathcal S}_{\ck}(0,A_k,A_{k-1})\; &\mbox{if } j<k,\\
 d_{\ck}^{\mathcal S}(J,A_k,A_{k-1})&\mbox{if } j=k>1.
 \end{cases}
\end{equation}
If $\max\{n,k\}=1$ or $j=k=1$, then $d_{\hi}^{\mathcal P}(P(\lambda))=\infty$.
\end{enumerate}
\end{theorem}
Before we give the proof we make a few remarks.
\begin{remark}{\rm
We observe the following simple facts about the inclusions
$\mathcal{P}^n_{k,j}\subseteq\mathcal{P}^n_{k+1,j}$, $k\geq 0$.
\begin{itemize}
\item[1)] It is an immediate corollary from equation \eqref{distdist1} that for $P(\lambda)\in\mathcal{P}^n_{k,j}$ one has
\[
d_{\sing}^{\mathcal{P}^n_{k,j}}(P(\lambda))= d_{\sing}^{\mathcal{P}^n_{\ell,j}}(P(\lambda)),\quad \ell\geq k.
\]
\item[2)] Observe that the closest singular polynomial may be of lower degree than the original one, e.g.,  let $\eps>0$ be small and let
\[
P(\lambda)=\lambda \mat{cc} \eps & 0\\ 0 & 0 \rix + \mat{cc} 0 & 0 \\ 0 & 1 \rix\in\mathcal{P}^{2}_{1,0}.
\]
Then the closest singular pencil in $\mathcal{P}^2_{1,0}$ is obtained by removing the $\eps$ entry, and this is a pencil of degree zero.
\item[3)] The (algebraic, geometric, partial) multiplicities of the eigenvalue infinity and the index of a matrix polynomial
$P(\lambda)\in\mathcal P_{k,j}^n$ are invariants with respect to the parameter $k$ and not with respect to the degree of the polynomial.
For example consider the matrix polynomial $P(\lambda)=\sum_{i=0}^k\lambda^iA_i\in\mathcal P^n_{k,j}$ with $A_0=I_n$ and $A_i=0$ for
$i=1,\dots,k$ which
is a matrix polynomial of degree zero. If $P(\lambda)$ is considered to be a matrix pencil (i.e. $k=1$) then it is of index one and the
algebraic multiplicity of $\infty$ is $n$. If, however, $P(\lambda)$ is considered as a quadratic matrix polynomial (i.e., $k=2$), then
it companion linearization has the form
\[
\lambda\mat{cc}0&0\\ 0&I_n\rix+\mat{cc}0&I_n\\ -I_n&0\rix
\]
and it follows that the eigenvalue $\infty$ has algebraic multiplicity $2n$ and the index is two.
The fact that a consistent spectral theory of matrix polynomials is only possible if the leading coefficients are allowed to be zero
is a well-known fact in the theory of matrix polynomials (see \cite{GohKL88}) and led to the introduction of the notion \emph{grade}
for the parameter $k$ in \cite{MacMMM11a}. Consequently, there is in general no equality between
$d_{\hi}^{\mathcal{P}^n_{k,j}}(P(\lambda))$ and $d_{\hi}^{\mathcal{P}^n_{\ell,j}}(P(\lambda))$ for $\ell>k$ which is also
reflected by formula \eqref{distdist2}.
\end{itemize}
}
\end{remark}
\proof
\eqref{tg1} The implication (a)$\Rightarrow$ (b) is trivial. Next, if $P(1)=-J+A_0+\cdots+A_k$ is singular,
then it follows from Proposition~\ref{charsinggen-new} that the kernels of $J,A_0,\dots,A_k$ have a nontrivial
intersection which, in turn, implies that $P(\lambda)$ is singular as obviously $\det P(\lambda)\equiv 0$.
Then \eqref{tg1a} is an easy consequence of \eqref{tg1}.

\eqref{tg2} First, consider the case $n=1$. If $k=1$, then $\mathcal I_{\hi}$ is clearly empty. If $k>1$
and $P(\lambda)\in\mathcal K$, then $J=A_k=A_{k-1}=0$ and the companion form of $P(\lambda)$ is
\[
L(\lambda)=\lambda\mat{cccc} 0&&&\\ &1&&\\ &&\ddots &\\ &&&1\rix+\mat{cccc}0&A_{k-2}&\dots&A_0\\ -1&0&&\\ &\ddots&\ddots&\\ &&-1&0\rix.
\]
By Lemma~\ref{lem:index} the index of $L(\lambda)$ and hence that of $P(\lambda)$ is at least two. If $P(\lambda)$ is regular, we thus
have $P(\lambda)\in\mathcal I_{\hi}$. If $P(\lambda)$ is singular, then it is identically zero and replacing $A_0$ with $\varepsilon$
and letting $\varepsilon\to 0$, we see that $P(\lambda)$ is in the closure of $\mathcal I_{\hi}$.

For $n>1$ we distinguish the cases $j<k$ and $j=k$.\\
\emph{Case $j<k$.} First observe that $\mathcal K$ is a closed set in $\mathcal P^n_{k,j}$ by Lemma \ref{closedJ}.
Hence, to prove the inclusion $\overline{\mathcal I^{\rule{0mm}{2mm}}
}_{\hi} \subseteq \mathcal K$ it suffices to show
that any matrix polynomial $P(\lambda)\in\mathcal I_{\hi}$ satisfies $\ker A_k\cap\ker A_{k-1}\neq\{0\}$.
To do this, suppose on the contrary that for some $P(\lambda)\in\mathcal I_{\hi}$ we have $\ker A_k\cap\ker A_{k-1}=\{0\}$.
If $\ker A_k=\{ 0\}$ then the matrix polynomial has no infinite eigenvalues and hence is of index zero. Hence, we may assume that
$A_k$ has a nontrivial kernel, and then there exists an orthogonal congruence transformation so that
\begin{equation}\label{pen1}
U^\top A_k U=\mat{cc} \widetilde A_k & 0 \\ 0 & 0\rix, \ U^\top J U=\mat{cc}J_{11} & J_{12} \\ -J_{12}^\top & J_{22} \rix,
\quad\mbox{and}\quad U^\top A_{k-1} U =\mat{cc}A_{11} & A_{12}  \\  A_{12}^\top & A_{22}  \rix,
\end{equation}
where $\widetilde A_k\in\Real^{n-r,n-r}$ ($r>0$) is invertible and all three matrices are partitioned conformably. In fact,
replacing $P(\lambda)$ by $U^\top P(\lambda) U$ if necessary, we may
assume $U=I_n$ in what follows. Since  $\ker A_k\cap\ker A_{k-1}=\{0\}$ and since $A_{k-1}$ is positive semidefinite,
it then follows that $A_{22}$ is invertible.

If $j<k-1$ then the companion linearization~\eqref{companion} of $P(\lambda)$ takes the form
\begin{equation}\label{compjsmall1}
\lambda\mat{cc|ccc} \widetilde A_k&&\\ &0&\\\hline &&I_{(k-1)n}\rix+
\mat{cc|c}A_{11}&A_{12}&\ast\\ A_{12}^\top&A_{22}&\ast\\\hline \ast &\ast&\ast\rix,
\end{equation}
where in comparison to~\eqref{companion} the first block row and column have been split into two, the last $k-1$ block rows and
columns have been merged into one, respectively, and $\ast$ denotes a possibly nonzero block entry.
Then it follows from Lemma~\ref{lem:index} (applied to the pencil that is obtained from~\eqref{compjsmall1} by permuting the
second and third block rows and columns) that the companion pencil and hence the matrix polynomial $P(\lambda)$ is of index one,
since $A_{22}$ is invertible.

If $j=k-1$, then the coefficient of $\lambda^{k-1}$ in $P(\lambda)$ is given by $A_{k-1}-J$ and hence the
companion linearization~\eqref{companion} of $P(\lambda)$ has the form
\begin{equation}\label{compjlarge}
\lambda\mat{cc|c} \widetilde A_k&&\\ &0&\\\hline &&I_{(k-1)n}\rix+
\mat{cc|c}A_{11}-J_{11}&A_{12}-J_{12}&\ast\\ A_{12}^\top+J_{12}^\top&A_{22}-J_{22}&\ast\\\hline \ast &\ast&\ast\rix
\end{equation}
with the same conventions as for the pencil~\eqref{compjsmall1}.
But with $A_{22}$ invertible also $A_{22}-J_{22}$ is invertible (see Proposition~\ref{charsinggen-new} applied with $\ell=0$ to $J_{22}$
and $X_0=A_{22}$). Again, it follows from Lemma~\ref{lem:index} that
the matrix polynomial $P(\lambda)$ is of index one.

For the converse inclusion $\mathcal K\subseteq\overline{\mathcal I^{\rule{0mm}{2mm}}}_{\hi}$ consider a matrix polynomial
$P(\lambda)\in\mathcal K$. Furthermore, let $u\in\ker A_k\cap\ker A_{k-1}$ with $\| u\|_2=1$ and let
$U\in \Real^{n,n}$ be orthogonal with last column $u$. Then
\begin{equation}\label{12.9.19}
U^\top A_k U=\mat{cc} \widetilde A_k & 0 \\ 0 & 0\rix,\ U^\top A_{k-1} U=\mat{cc} \widetilde A_{k-1} & 0 \\ 0 & 0\rix, \ U^\top J U
=\mat{cc} J_{11} & v \\ -v^\top & 0\rix,
\end{equation}
with $\widetilde A_k,\widetilde A_{k-1},J_{11}\in\mathbb R^{n-1,n-1}$ (not necessarily being invertible) and $v\in\mathbb R^{n-1}$.
Note that the entry $0$ in $(2,2)$ block of $ U^\top J U$ is caused by the skew-symmetry of $J$. Again,  replacing
$P(\lambda)$ with $U^\top P(\lambda)U$ if necessary, we may assume without loss of generality that $U=I_n$.

First assume that $j<k-1$. For small $\eps>0$ we have that  $\widetilde A_k+\eps I_{n-1}\in\mathbb R^{n-1,n-1}$ and
$A_j-J+\eps I_n\in\mathbb R^{n,n}$ are invertible.  Then the matrix polynomial
$P_\eps(\lambda)$ that is obtained from $P(\lambda)$ by replacing $A_k$ with $A_k+\diag(\eps I_{n-1},0)$ and
$A_j-J$ with $A_j-J+\eps I_n\in\mathbb R^{n,n}$ is regular, because at least one coefficient (namely the coefficient
associated with $\lambda^j$) is invertible. Furthermore, the companion linearization of $P_\eps(\lambda)$ takes the form
\begin{equation}\label{compjsmall0}
\lambda\mat{cc|c} \widetilde A_k+\eps I_{n-1}&&\\ &0&\\\hline &&I_{(k-1)n}\rix+
\mat{cc|c}\widetilde A_{k-1}&0&\ast\\ 0&0&\ast\\\hline \ast &\ast&\ast\rix.
\end{equation}
By Lemma ~\ref{lem:index} we see that this pencil, and hence $P_\eps(\lambda)$ itself, has  index greater than one.

Now let $j=k-1$. For sufficiently small $\eps>0$
we have that $\widetilde A_k+\eps I_{n-1}\in\mathbb R^{n-1,n-1}$ and $A_{k-1}+\eps I_n\in\mathbb R^{n,n}$ are invertible
and that $v+\eps e_1\neq 0$. Let $P_\eps(\lambda)$ be the matrix polynomial obtained from $P(\lambda)$ by replacing $A_k$ and with
$A_k+\diag(\eps I_{n-1},0)$ and $A_{k-1}$ with $A_{k-1}+\eps I_n$ as well as
$v$ in $J$ with $v+\eps e_1\neq 0$. Again, it follows that $P_\eps(\lambda)$ is regular as the coefficient at $\lambda^{k-1}$ is
regular (see Proposition~\ref{charsinggen-new} applied with $\ell=0$ to $J$ and $X_0=A_{k-1}+\eps I_{n}$).
Then the companion linearization of $P_\eps(\lambda)$ takes the form
\begin{equation}\label{compjsmall}
\lambda\mat{cc|c} \widetilde A_k+\eps I_{n-1}&&\\ &0&\\\hline &&I_{(k-1)n}\rix+
\mat{cc|c}\widetilde A_{k-1} +\eps I_{n-1}- J_{11}&-v-\eps e_1&\ast\\ v^\top+\eps e_1^\top&0&\ast\\\hline \ast &\ast&\ast\rix.
\end{equation}
By Lemma ~\ref{lem:index} we see that~\eqref{compjsmall}, and hence $P_\eps(\lambda)$ itself, has  index greater than one.

Letting $\eps\to 0$ we see  that in both cases $j<k-1$ and $j=k-1$ we have $\mathcal I_{\hi}\ni P_{\eps}(\lambda)\to P(\lambda)\in
\overline{\mathcal I^{\rule{0mm}{2mm}}}_{\hi}$.


\emph{Case $j=k$.}
If $j=k=1$ and $P(\lambda)=\lambda(J+A_1)+A_0\in\mathcal P^n_{1,1}$, then by Theorem~\ref{thm:singind}, zero is a semisimple eigenvalue
(if it is an eigenvalue) of the reversal $\lambda A_0+J+A_1$ of $P(\lambda)$ and hence the index of $P(\lambda)$ is at most one.
This shows that $\mathcal I_{\hi}$ is empty in that case.

Thus, assume that $j=k>1$. To show the inclusion $\overline{\mathcal I^{\rule{0mm}{2mm}}}_{\hi} \subseteq \mathcal K$
suppose, as in the proof of \eqref{tg2},  that for some $P(\lambda)\in\mathcal I_{\hi}$ we have $\ker J\cap\ker A_k\cap\ker A_{k-1}=\{0\}$.
If $\ker (A_k-J)=\{ 0\}$ then the matrix polynomial has no infinite eigenvalues and hence is of index zero. Hence, we may assume that
$A_k-J$ has a nontrivial kernel, which by Proposition~\ref{charsinggen-new} applied with $\ell=0$ to $J$ and $X_0=A_k$ implies that
 there exists an orthogonal congruence transformation $U$ so that
\begin{equation}\label{pen2}
U^\top A_k U=\mat{cc} \widetilde A_k & 0 \\ 0 & 0\rix, \ U^\top J U=\mat{cc}J_{11} & 0 \\ 0 &0 \rix,
\quad\mbox{and}\quad U^\top A_{k-1} U =\mat{cc}A_{11} & A_{12}  \\  A_{12}^\top & A_{22}  \rix,
\end{equation}
where $\widetilde A_k-J_{11}\in\Real^{n-r,n-r}$ ($r>0$) is invertible and all three matrices are partitioned conformably. In fact,
replacing $P(\lambda)$ by $U^\top P(\lambda) U$ if necessary, we may
assume $U=I_n$ in what follows. Since  $\ker J\cap \ker A_k\cap\ker A_{k-1}=\{0\}$ and since $A_{k-1}$ is positive semidefinite,
it follows that $A_{22}$ is invertible.
The companion matrix pencil \eqref{companion} of $P(\lambda)$ takes the form
\begin{equation}\label{compjsmall2}
\lambda\mat{cc|ccc} \widetilde A_k-J_{11}&&\\ &0&\\\hline &&I_{(k-1)n}\rix+
\mat{cc|c}A_{11}&A_{12}&\ast\\ A_{12}^\top&A_{22}&\ast\\\hline \ast &\ast&\ast\rix
\end{equation}
and from Lemma~\ref{lem:index} we infer that
the matrix polynomial $P(\lambda)$ is of index one.
%
%

For the converse $\mathcal K\subseteq\overline{\mathcal I^{\rule{0mm}{2mm}}}_{\hi}$ let $P(\lambda)\in\mathcal K$. Then
we have the following decomposition
\[
U^\top A_k U=\mat{cc} \widetilde A_k & 0 \\ 0 & 0\rix,\ U^\top A_{k-1} U=\mat{cc} \widetilde A_{k-1} & 0 \\ 0 & 0\rix,
\ U^\top J U=\mat{cc} J_{11} & 0 \\ 0 & 0\rix,
\]
\[
\mbox{and }\ U^\top A_{k-2}U=\mat{cc}  A_{11} & A_{12} \\ A_{12}^\top & A_{22}\rix, \ U^\top J U=\mat{cc} J_{11} & 0 \\ 0 & 0\rix,
\]
with $A_k,A_{k-1},A_{11},J_{11}\in\Real^{n-1,n-1}$ not necessarily invertible.
Replacing $\widetilde A_{k-2}$ with $\widetilde A_{k-2}+\eps I_{n-1}$, then for sufficiently small $\varepsilon$ we
we get a family of regular pencils $P_\eps(\lambda)$ of index at least two and
such that $\mathcal I_{\hi}\ni P_{\eps}(\lambda)\to P(\lambda)\in \overline{\mathcal I^{\rule{0mm}{2mm}}}_{\hi}$.

\eqref{tg2a} is an immediate consequence of \eqref{tg2}.
\eproof

\begin{remark}{\rm
At first, it may come as a surprise that a matrix polynomial $P(\lambda)$ as in Theorem~\ref{distances} is already
singular if $P(1)$ is singular which means that $1$ cannot be an eigenvalue of a regular $P(\lambda)\in\mathcal P^n_{k,j}$.
More generally, if $\alpha>0$ then replacing $\lambda$ with
$\frac{\lambda}{\alpha}$ and $A_i$ with $\alpha^iA_i$ shows that if $P(\alpha)$ is singular, then $P(\lambda)$ is already
a singular matrix polynomial. This generalizes in a nontrivial way the observation that any scalar polynomial with
nonnegative coefficients cannot have real zeros that are positive unless it is the zero polynomial.}
\end{remark}

\begin{remark}\label{rem_8.11.19}{\rm
The reason for not investigating the structured  distance to instability for $P(\lambda)$  in Theorem~\ref{distances} is
the fact that in contrast to  Theorem ~\ref{singind}\eqref{singindII} the distances to higher index and
instability need not coincide for matrix polynomials of degree larger than one. We will return to the distance to instability for
quadratic polynomials in Section~\ref{sec:qp}, because that task
is still accessible by the common kernel methods framework. This is due to a nontrivial result, Theorem \ref{thm:24.10.19} below,
which states that the only spectral points that may cause instability are zero and infinity. However, already for  degree three
the reason for instability may be different, since a polynomial in $\mathcal{P}^n_{3,\ell}$ might have eigenvalues in the right half plane.
For example, the scalar polynomial $p(\lambda)=\lambda^3+1\in\mathcal P^1_{3,\ell}$
 (thus $J=X_1=X_2=0$, $X_0=X_3=1$) has  two  roots  in the right half plane.}
\end{remark}

\section{Distance problems for first order dHDAE systems}\label{sec:dhdaedist}

In this section we will revisit the distance problems for first order dHDAE systems formulated in Section~\ref{sec:prob}.
We will first present the missing proofs which are now easy consequences of the extended results from the previous
section. Then, we will present two examples that show that the structured distances for dHDAE
systems may differ considerably from the corresponding ones under arbitrary perturbations.

\subsection{Proofs of and comments on the main results in Section~\ref{sec:prob}}\label{s:proof}

\noindent
{\bf Proof of Theorem \ref{Q=I}.}
\eqref{I4} It  follows from Theorem \ref{distances}\eqref{tg1} ($k=1$, $j=0$) that $P(\lambda)$ is singular if and only if the kernels
of $E,J$ and $R$ have a nontrivial intersection. Applying Theorem
\ref{charsinggen-new} ($\ell=1$) we get the desired transformation $U$.

\eqref{I3} By Theorem~\ref{thm:singind} the index is at most two.
Then it  follows from Theorem \ref{distances}\eqref{tg2} ($k=1$, $j=0$) that $P(\lambda)$ is in the closure of regular dH pencils of
index 2  if and only if the kernels of $E$ and $R$ have a nontrivial intersection.
\eproof

\medskip

\noindent
{\bf Proof of Theorem \ref{singind}.}
\eqref{singindI} The proof is obtained from Theorem~\ref{main} with $k=1$, $X_0=E$ and $X_1=R$ and Theorem
 \ref{distances}, and using $J^\top J=-J^2$.

\eqref{singindII} First, it immediately follows from Theorem~\ref{thm:singind} that $P(\lambda)$ is stable if and only if it
is regular and of index one.
The proof is then obtained from Theorem~\ref{main} with $k=1$, $X_0=E$ and $X_1=R$ and using $J^\top J=-J^2$.
By Theorem~\ref{thm:singind} any pencil  $\lambda E-(J-R)$ with $E,R\geq 0$ and $J^\top=-J$, associated with a dH system,
is of index at most $2$.
\eproof

\medskip

We have the following immediate corollary of Theorems \ref{Q=I} and \ref{singind}.
 \begin{corollary}\label{closures}
 For  $J=-J^\top$, $E,R\geq 0$ one has the following estimate
\[
2\cdot \lambda_{\min}(-J^2) +d_{\hi}^{\mathcal{L}}\big(\lambda E - (J-R)\big)^2\leq d_{\sing}^{\mathcal{L}}\big(\lambda E - (J-R)\big)^2\leq 2
\norm{J}^2+d_{\hi}^{\mathcal{L}}\big(\lambda E - (J-R)\big)^2.
\]
In particular, the set of singular dH pencils in $\mathcal{P}^n_{1,0}$ is contained in the closure
in $\mathcal P^n_{1,0}$ of the set of index two  regular dH pencils.
\end{corollary}

\begin{remark}\label{permuting}{\rm Consider the reversed pencil $-\lambda (J-R)+E$ with $J=-J^\top$ and $E,R\geq 0$.
Statement \eqref{tg2} of Theorem~\ref{distances} shows that
its structured distance to higher index $d^{\mathcal{L}}_{\hi}\big(-\lambda (J-R)+E \big)$ equals its distance to singularity
$d_{\sing}^{\mathcal{L}}\big(-\lambda (J-R)+E \big)$.    This is in line with the fact that $\infty$ can only be a semisimple
eigenvalue of $-\lambda (J-R)+E$
as zero is a semisimple eigenvalue of  $\lambda E-(J-R)$, cf. Theorem \ref{Q=I}. Then, in summary, we have
\begin{align*}
d^{\mathcal{P}}_{\ck}( J,E,R  )&=
d^{\mathcal{P}}_{\hi}\big(-\lambda (J-R)+E \big)\\
&=d^{\mathcal{P}}_{\sing}\big(-\lambda (J-R)+E \big)\\
&=d^{\mathcal{P}}_{\sing}\big(\lambda E- (J-R) \big)\\
&=d^{\mathcal{P}}_{\sing}\big(\lambda R - (J-E) \big).
\end{align*}
}
\end{remark}
Since our main focus is on distance problems in this paper, it was necessary to characterize the
closure of the set of dH pencils of index two in Theorem \ref{Q=I}~\eqref{l3}. Our next result gives a characterization of the set of
regular dH pencils of index two.

\begin{proposition}  A pencil $L(\lambda)=\lambda E-(J-R)$ with $E,R\geq 0$, $J=-J^\top$ is regular of index two if and only if there
exists an orthogonal matrix  $U$ such that
\begin{equation}\label{decoi2reg}
U^\top EU=\mat{ccc} E_{11} & 0 & 0 \\ 0 & 0 & 0 \\ 0 & 0 & 0\rix,\quad U^\top JU=\mat{ccc} J_{11} & J_{12} & J_{13} \\ -J_{12^\top} & J_{22}
& 0 \\ -J_{13}^\top & 0 & 0\rix,\quad U^\top RU=\mat{ccc} R_{11} & R_{12} & 0 \\ R_{12}^\top & R_{22} & 0 \\ 0 & 0 & 0\rix,
\end{equation}
where $n=p+q+r$, $p,r> 0$, $E_{11}\in\Real^{p,p}$ is invertible, $J_{22}-R_{22}\in\Real^{q,q}$ is invertible, and $J_{13}\in\Real^{p,r}$
has full column rank.
\end{proposition}
\proof
Suppose that the decomposition \eqref{decoi2reg} holds. As $J_{22}-R_{22}$ is invertible and $J_{13}$ has full column rank, we have that
$-J+E+R$ is invertible. Hence, by Theorem \ref{Q=I}\eqref{I4} and Theorem \ref{distances}\eqref{tg1} the pencil is regular.  Therefore,
by Lemma \ref{lem:index} it is of index two. To prove the converse implication first we find an orthogonal transformation $U_1$
such that
 \begin{equation}\label{decoi3reg}
U^\top _1EU_1=\mat{ccc} E_{11} & 0  \\ 0 & 0  \rix,\quad U^\top _1JU_1=\mat{ccc} J_{11} & J' \\ -J'{}^\top & \widetilde J \rix,
\quad U^\top _1RU_1=\mat{ccc} R_{11} & R'  \\ R'{}^\top & \widetilde R \rix,
 \end{equation}
with $E_{11}\in\Real^{p,p}$ invertible. By Lemma \ref{lem:index} the matrix $\widetilde J-\widetilde R$ is singular. Applying
Theorem \ref{charsinggen-new} to $\widetilde J$ and $\widetilde R$ we get an orthogonal transformation and
a splitting of the last $n-p$ rows and columns, which combined with \eqref{decoi3reg} gives
\begin{equation}\label{decoi4reg}
U^\top EU=\mat{ccc} E_{11} & 0 & 0 \\ 0 & 0 & 0 \\ 0 & 0 & 0\rix,\; U^\top JU=\mat{ccc} J_{11} & J_{12} & J_{13} \\ -J_{12^\top} & J_{22}
& 0 \\ -J_{13}^\top & 0 & 0\rix,\; U^\top RU=\mat{ccc} R_{11} & R_{12} & R_{13} \\ R_{12}^\top & R_{22} & 0 \\ R_{13}^\top & 0 & 0\rix,
\end{equation}
with some orthogonal $U$. As $R$ is positive semidefinite we have $R_{13}=0$. But then $J_{13}$ needs to have full column rank, otherwise
the pencil would be singular.
\eproof
\subsection{Structured vs. unstructured distances}\label{sec:strunstr}
In this subsection we compare the structured and the unstructured distances.
First note that the statements of Theorem~\ref{Q=I} are not true without the structure assumptions on $E$ and $R$.
While it is obvious that~\eqref{I4} cannot hold for arbitrary pencils, we need a simple example to disprove an unstructured analogue
of~\eqref{I3}.
\begin{example}\label{ex2}{\rm
Let $n=2$, $E=J=0$, $R=I_2$. Then setting $\widetilde J=J$, $\widetilde R=R$ and
$$
\widetilde E= \mat{cc} 0 & \eps \\ 0 & 0 \rix
$$
we see that $\lambda \widetilde E-(\widetilde J-\widetilde R)$ is regular and of index two (but not dissipative Hamiltonian).
Letting $\eps\to 0$, we find that $\lambda E-(J-R)$ is in the closure of regular pencils of index two
although $\ker{E}\cap\ker{R}=\{0\}$.}
\end{example}

To analyze the distances in Theorem~\ref{singind}, recall that in \cite{ByeHM98} the (unstructured) distance to singularity was defined as
\[
d_{\sing}(\lambda E -A):=\inf\big\{\|[\Delta_E,\Delta_A]\|_F\ \big| \lambda(E+\Delta_E)+A+\Delta_A \text{ is singular}\big\}.
\]
\begin{example}\label{numeric1}{\rm
Let
\[
E = \mat{cc}  0 &    0\\
     0   &  1\rix,\quad
J = \mat{cc} 0 &  -0.5 \\
    0.5  &    0 \rix,
\quad
R = \mat{cc}  0.18 &   0.42\\
    0.42 &   1.03\rix\geq 0.
\]
Then (rounding the numerical results to four digits) we have
\[
\lambda_{\min}(-J^2+E^2+R^2)=0.5819,\quad \sigma_{\min}\left(\mat{c} A\\E\rix\right)=0.1908,\quad
\sigma_{\min}\left(\mat{cc} A&E\rix\right)=0.6056,
\]
where $\sigma_{\min}$ stands for the smallest singular value.  The first equality implies, by Theorem \ref{singind}\eqref{singindII},
that $d_{\sing}^{\mathcal{L}}(\lambda E-(J-R)))\geq 0.5819$, while the second and third equality imply, together with Corollary 3
of \cite{ByeHM98}, that $d_{\sing}(\lambda E-A)=0.1908$.
}
\end{example}

The next example shows mainly the same behaviour, though we refine the constraints. Namely,
we show that if we change the constraint in the definition of \eqref{distEJR} to
 $E+\Delta_ E,R+\Delta_R$ being symmetric (instead of being positive definite)
then we get an essentially different distance.
\begin{example}\label{ex1}{\rm
Consider a dissipative Hamiltonian system \eqref{dHwoQ} with coefficients
\[
	E=\matp{1\\&1\\&&1\\&&&0\\&&&&0},\quad
	J=\matp{0&0&0&0&-1\\ 0&0&0&1&1\\ 0&0&0&-1&-1\\0&-1&1&0&\eps\\1&-1&1&-\eps & 0},\quad
R=\matp{\alpha&0&0&0&1\\0&\alpha&0&1&1\\0&0&\alpha&1&1\\0&1&1&\eps &0\\1&1&1&0&\eps},
\]
where $\eps>0$ and
\[
	\alpha=\eps^{-1}\cdot\norm{\matp{0&1\\1&1\\1&1}}^2_F+1.
\]
By Theorem 1.1 of \cite{FoiF90} such choice of $\alpha$ makes $R>0$. Consider now the perturbation
\[
  \Delta_E=0,\quad	\Delta_J=\matp{0_{3\times 3} & \\&& -\eps\\ &\eps},\quad
	 \Delta_R=\matp{0_{3\times 3} & \\& -\eps\\ &&-\eps}.	
\]
Clearly $\Delta_J=-\Delta_J^\top$ and $\Delta_R=\Delta_R^\top$, and $\Delta_E=\Delta_E^\top$. The pencil $\lambda E-\widehat A$ with
\[
\widehat A:=J+\Delta_J-(R+\Delta_R)=\matp{-\alpha&0&0&0&-2\\0&-\alpha&0&0&0\\0&0&-\alpha&-2&-2\\
0&-2&0&0& 0\\0&-2&0&0&0}
\]
is now singular, because one easily checks that $\det(\lambda E-\widehat A)\equiv 0$.
In this way, we have constructed a perturbation with
\[
		\norm{ [ \Delta_J,\Delta_R,\Delta_E] }_F=2\eps,
\]
such that the perturbed pencil is singular, but the perturbation is not structure-preserving, because
the matrix $R+\Delta_R$ is now indefinite.  Observe also that, unlike in Theorem \ref{Q=I}\eqref{I4},
$E$ and $\widehat A$ do not have a common right or left kernel. This means that in the Kronecker canonical form there are left and
right minimal indices of size at least one.
		
On the other hand, we have
\[
J^\top J+E^2+R^2=\matp{3+\alpha^2&0&2&-\eps&\alpha+\eps\\ 0&5+\alpha^2&0&\alpha+2\eps&\alpha\\
2&0&5+\alpha^2&\alpha&\alpha+2\eps\\ -\eps&\alpha+2\eps&\alpha&4+2\eps^2&4\\ \alpha+\eps&\alpha&\alpha+2\eps&4&6+2\eps^2}
\]
and a simple {\sc Matlab} calculation shows that we have
$\lambda_{\min}(J^\top J +E^2+R^2)^{1/2}\geq 0.81$ for $\eps\in(10^{-1},10^{-6})$.
Hence, the smallest perturbation $\Delta_E,\Delta_J,\Delta_R$ that makes the pencil singular, while keeping
$R+\Delta_R\geq0$, $E+\Delta_E\geq0$, and $J+\Delta_J$ skew-symmetric, satisfies
\[
	\norm{ [ \Delta_J,\Delta_R,\Delta_E] }_F\geq 0.81,
\]
for $\eps\in(10^{-1},10^{-6})$ by Theorem \ref{singind}\eqref{singindI}.
}
\end{example}

Note that the pencil $\lambda E-(J-R)$ from Example \ref{ex2} also shows
that Theorem~\ref{singind}\eqref{singindII} does not hold for
unstructured perturbations. Indeed, that pencil is in the closure of regular pencils of index 2, but
$d_{\hi}^{\mathcal{P}}( \lambda E-(J-R)) \geq \lambda_{\min} (E^2+R^2)=1$.

As last example we consider the analysis of distance problems in circuit simulation.
\begin{example}\label{exRCL}{\rm
A simple RLC network, see, e.g., \cite{BeaMXZ18,Dai89,Fre08,Fre11}, can be modeled by a dHDAE system of the
form
\begin{eqnarray}\label{eq:RLC_network_1}
			\underbrace{\mat{ccc} G_c C G_c^\top & 0 & 0\\ 0 & L & 0 \\ 0 & 0 & 0 \rix}_{=:E}
			\mat{c}\dot{V} \\ \dot{I_l}\\ \dot{I_v} \rix=
			\underbrace{\mat{ccc} -G_r R_r^{-1}G_r^\top & -G_l & -G_v \\ G_l^\top & 0 & 0\\ G_v^\top & 0 & 0 \rix}_{=:J-R}
			\mat{c}V \\ I_l\\ I_v \rix,
\end{eqnarray}
where $L > 0$, $C > 0$, $R_r> 0$ are real symmetric matrices describing inductances, capacitances, and resistances, respectively.
The subscripts $r,\,c,\,l,$ and $v$  refer  to the resistors, capacitors, inductors, and voltage sources,
while $V$, $I$ denote voltage and current, respectively. The matrices $G_c,G_l,G_r,G_v$ encode the network topology, see~\cite{Fre08} for
details. Here, $J$ and $-R$ are defined to be the skew-symmetric and symmetric parts, respectively, of the matrix on the right hand side
of~\eqref{eq:RLC_network_1}. We see that $E,R\geq 0$ and $J=-J^T$.
			
It was shown in \cite[Theorem 1]{Fre08} that the pencil $\lambda E-(J-R)$ is regular if and only if $G_v$ has full column rank and
$$
G_1:=\matp{  G_c & G_r & G_l & G_v}
$$
has full row rank. Note that this equivalence is now a simple corollary of Proposition~\ref{singind}\eqref{singindI}. Indeed, $\lambda E-(J-R)$
is singular if and only if the kernels of the three matrices
\[
E=\mat{ccc} G_cCG_c^\top&0&0\\  0&L&0\\ 0&0&0\rix,\quad J=\mat{ccc}0&-G_l&-G_v\\ G_l^\top &0&0\\ G_v^\top &0&0\rix,\quad\mbox{and}\quad
R=\mat{ccc} G_rR_r^{-1}G_r^\top &0&0\\ 0&0&0\\ 0&0&0\rix
\]
have a nontrivial intersection. Having in mind that $C$, $L$, and $R_r^{-1}$ are positive definite matrices, we immediately obtain
that $x=\mat{ccc}x_1^\top&x_2^\top&x_3^\top\rix^\top\in\ker{E}\cap\ker{J}\cap\ker{R}$ if and only if
$G_c^\top x_1=0$, $x_2=0$ as well as $G_l^\top x_1=0$, $G_v^\top x_1=0$, $G_v x_3=0$, and $G_r^\top x_1=0$ which in turn is equivalent to
\[
x_1^\top G_1=0,\quad x_2=0,\quad\mbox{and}\quad G_v x_3=0.
\]
Thus, we see that $\lambda E-(J-R)$ is singular if and only if either $G_1$ does not have full row rank or $G_v$ does not have
full column rank.

All this shows that regularity of the pencil $\lambda E-(J-R)$ depends only on the network topology, cf. Remark 1 in~\cite{Fre08}.
As one can expect,
the distance to singularity depends also on the values of matrices $L,C,R_r$. Observing that the matrix $-J^2+R^2+E^2$ has
the form
\[
-J^2+R^2+E^2=\mat{ccc}(G_cCG_c^\top)^2+(G_rR_r^{-1}G_r^\top)^2+G_lG_l^\top+G_vG_v^\top &0&0\\ 0&L^2+G_l^\top G_l & G_l^\top G_v\\
0&G_v^\top G_l&G_v^\top G_v\rix,
\]
we obtain by~\eqref{lmino} that the structured distance to singularity  is bounded by
\[
\lambda_{\min}(-J^2+R^2+E^2)^{1/2}\leq d_{\sing}^{\mathcal{L}}\big(\lambda E - (J-R)\big)\leq 2\cdot\lambda_{\min}(-J^2+R^2+E^2)^{1/2},
\]
where $\lambda_{\min}(-J^2+R^2+E^2)$ is the minimum of the two values
\[
\lambda_{\min}\big( (G_cCG_c^\top)^2 +(G_rR_r^{-1} G_r^\top)^2+G_lG_l^\top+G_vG_v^\top\big)\quad\mbox{and}\quad
\lambda_{\min}\left(\mat{cc}L^2+G_l^\top G_l & G_l^\top G_v\\ G_v^\top G_l&G_v^\top G_v\rix\right).
\]
}
\end{example}

In this section we have characterized the distances to singularity, high index and instability for linear constant coefficient dH systems.
In the next section we extend these results to quadratic matrix polynomials.

\section{Quadratic polynomials, linearization, and removing $Q$}\label{sec:qdh}

In this section, we will show how the main results of the previous sections can be applied to more general
situations. We first study how the transformation (linearization) of structured matrix polynomials to dH pencils can be performed and
then apply the distance results to quadratic matrix polynomials. Finally we discuss the more general dH pencils in \eqref{dH} and
show how the multiplier $Q$ can be removed.

\subsection{Dissipative Hamiltonian linearisations}\label{sec:dhlin}
Consider a second order system of the form
\begin{equation}\label{second}
M\ddot x-(G-D)\dot x+Kx=0
\end{equation}
with $M,G,D,K\in\mathbb R^{n,n}$ satisfying $M,D,K\geq0$ and $G=-G^\top$.
A companion linearization \eqref{companion} of \eqref{second} is given by
\[
L(\lambda)=\lambda\mat{cc}M&0\\ 0&I\rix-\mat{cc}D-G&K\\ -I&0\rix=\lambda\mat{cc}M&0\\ 0&I\rix-\left(\mat{cc}G&I\\ -I&0\rix
-\mat{cc}D&0\\ 0&0\rix\right)\mat{cc}I&0\\ 0&K\rix.
\]
It satisfies the hypothesis of Theorem~\ref{thm:singind} with
\[
E=\mat{cc}M&0\\ 0&I\rix,\quad J=\mat{cc}G&I\\ -I&0\rix,\quad R=\mat{cc}D&0\\ 0&0\rix,\quad Q=\mat{cc}I&0\\ 0&K\rix,
\]
and thus we immediately have the following result.
\begin{theorem}\label{thm:24.10.19}
Let $P(\lambda)=\lambda^2 M-\lambda (G-D)+K$, where $M,D,K,G\in\mathbb C^{n,n}$ with $G^H=-G$ and $M,D,K\geq 0$.
Then the following statements hold.
\begin{enumerate}
\item[\rm (i)] All eigenvalues of $P(\lambda)$ are in the closed left half complex plane and all finite nonzero
eigenvalues on the imaginary axis are semisimple.
\item[\rm (ii)] The possible length of Jordan chains of $P(\lambda)$, associated with either  the eigenvalue $\infty$ or the
eigenvalue zero, is at most two.
\item[\rm (iii)] All left and all right minimal indices of $P(\lambda)$ are zero (if there are any).
\end{enumerate}
\end{theorem}
\proof
The proof of (i) and the statement in (ii) on the length of the Jordan chains associated with the eigenvalue
$\infty$ follows immediately from Theorem~\ref{thm:singind}, while the remaining statement of (ii) then follows by applying the already
proved part of (ii)
to the reversal $\lambda^2K+\lambda (G-D)+M$ of $P(\lambda)$, which has the same structure.
To see (iii), observe that by Theorem~\ref{thm:singind} the left minimal indices of $L(\lambda)$
are all zero and the right minimal indices of $L(\lambda)$ are at most one. By \cite[Theorem 5.10]{DeDM09a}
the left minimal indices of $P(\lambda)$ coincide with those of $L(\lambda)$, and if $\ve_1,\dots,\ve_k$
are the right minimal indices of $P(\lambda)$, then $\ve_1+1,\dots,\ve_k+1$ are the right minimal indices
of $L(\lambda)$. This implies that all minimal indices of $P(\lambda)$ are zero (if there are any).
\eproof
\begin{remark}{\rm One may be tempted  to use the results on dH pencils with $Q=I$ to prove Theorem~\ref{thm:24.10.19}, i.e., to apply
Theorem \ref{Q=I} instead of Theorem~\ref{thm:singind} and also to extend the result to polynomials of higher degree.
However, as the constant term of the companion form \eqref{companion} contains a principal submatrix
\[
\mat{cc} Y_k & Y_{k-1} \\ -I_n & 0 \rix,
\]
 which has a positive semidefinite symmetric part only when $Y_{k-1}=I$, in general one needs to consider different linearizations,
 e.g., $L(\lambda)$ to be from
one of the classical linearizations classes in \cite{MacMMM06b,MacMMM06a} or  a so-called \emph{Fiedler linearization},
see, e.g., \cite{DeDM10a}.
However, one still meets the following general obstacles.
\begin{enumerate}[\rm 1)]
\item As we have seen in Remark~\ref{rem_8.11.19}, a polynomial from $\mathcal P^n_{k,j}$ may have spectrum in the right
half plane if $k>2$. Hence, no linearization of such a polynomial satisfies the hypothesis of Theorem~\ref{thm:singind}.
\item The index of a polynomial from $\mathcal P^n_{k,j}$  may be larger than two if $k>2$; consider e.g.
$P(\lambda)=\lambda^3 X_3+\lambda^2 X_2+\lambda X_1+\lambda X_0-J$ with $X_3=X_2=X_1=J=0$ and $X_0=1$. The
companion form is
\[
\lambda\mat{ccc}0&&\\ &1&\\ &&1\rix+\mat{ccc}0&0&1\\ -1&&\\ &-1&\rix,
\]
which corresponds to a block of size three at $\infty$, i.e., $P(\lambda)$ is of index 3. Hence, none of its index preserving (strong)
linearizations satisfies the hypothesis of Theorem~\ref{thm:singind}.
\item Let $P(\lambda)\in\mathcal P^n_{k,j}$ be a singular matrix polynomial with left minimal indices $\eta_1,\dots,\eta_\ell$ ($\ell>0$) and
right minimal indices $\varepsilon_1,\dots,\varepsilon_\ell$; note that the numbers of left and right minimal indices coincide, because the
matrix polynomial is square. Take  $L(\lambda)$ as any of the structured  linearizations in \cite{DeDM10a,MacMMM06b,MacMMM06a}. Then, there
exists a number $q\in\{0,1,2,\dots,k-1\}$, known in advance for a particular linearization, such that $L(\lambda)$ has the left minimal indices
$\eta_1+q,\dots,\eta_\ell+q$ and right minimal indices $\varepsilon_1+k-1-q,\dots,\varepsilon_\ell+k-1-q$, see \cite{DeDM09a,DeDM10a}.
Thus, even for $k=2$, the linearization $L(\lambda)$ cannot satisfy the hypothesis of Theorem \ref{Q=I} and for
$k>2$ it cannot satisfy the hypothesis of  Theorem~\ref{thm:singind}.
\end{enumerate}

Hence, if $k>2$, then a linearization $L(\lambda)$ of $P(\lambda)$ cannot be a dH pencil with arbitrary $Q$  and if  $k=2$ it can
only be a dH pencil with (necessarily nontrivial) $Q$. It remains an open question to derive additional conditions on the matrix coefficients
for a polynomial $P(\lambda)\in\mathcal P^n_{k,j}$
in order to guarantee that its spectrum is in the closed left half plane such that all eigenvalues on the imaginary axis
(possible excluding zero and infinity) are semisimple.
}
\end{remark}
In the next subsection we study the case of quadratic matrix polynomials with dH structure.

\subsection{Quadratic matrix polynomials with dH structure}\label{sec:qp}

In the case of quadratic matrix polynomials with dH structure, i.e., if $P(\lambda)=\lambda^2A_2+\lambda A_1+A_0 \in \mathcal P^n_{2,1}$,
where $A_2,A_0,A_1+A_1^\top\geq 0$, the structured distances to singularity, to higher index, and to instability were defined
in~\eqref{distP},~\eqref{indexdistP}, and~\eqref{instdistP}, respectively. Recalling that for a matrix $Y\in\mathbb R^{n,n}$ and a
vector $u\in\mathbb R^n$ with $\|u\|_2=1$
the perturbation matrix $\Delta_Y^u$ in~\eqref{derpert} is defined as  $\Delta_Y^{u}=-uu^\top Y-Yuu^\top+uu^\top Yuu^\top$,
we obtain the following result.
\begin{theorem}\label{thm:quadpoly}
Let $P(\lambda)=\lambda^2 M-\lambda (G-D)+K\in\mathcal{P}^n_{2,1}$, i.e., $M,D,K,G\in\mathbb R^{n,n}$ with $G^\top=-G$ and $M,D,K\geq 0$.
\begin{enumerate}[{\rm (i)}]
\item\label{pp1} The matrix polynomial $P(\lambda)$ is singular if and only if all four matrices $M,D,K,G$ have a common kernel.
In particular,
the structured distance to singularity $d_{\sing}^{\mathcal{P}}(P(\lambda))$ is attained for a perturbation of the form $\Delta_M^u$,
$\Delta_D^u$,
$\Delta_K^u$, $\Delta_G^u$ for some $u\in\mathbb R^{n,n}$ with $\|u\|_2=1$ and satisfies
\begin{equation*}
d_{\sing}^{\mathcal{P}}(P(\lambda))=d_{\ck}^{\mathcal{S}}(G,M,D,K),
\end{equation*}
where $d_{\ker}^{\mathcal{S}}(G,M,D,K)$ is given by Theorem~\mbox{\rm \ref{main}}. In particular, the structured distance to singularity is
bounded by
\begin{equation}\label{polydist}
\sqrt{\lambda_{\min} (M^2+D^2+K^2-G^2) }\leq d_{\sing}^{\mathcal{P}}(P(\lambda)) \leq \sqrt{2\cdot\lambda_{\min} (M^2+D^2+K^2-G^2)}.
\end{equation}
\item\label{pp2} The matrix polynomial $P(\lambda)$ is in the closure of polynomials in ${\mathcal{P}^n_{2,1}}$ that are regular and of
index larger than one (and thus of index exactly two)
if and only if $M$ and $D$ have a common kernel. In particular,
the structured distance to higher index $d_{\hi}^{\mathcal{P}}(P(\lambda))$ is attained for a perturbation of the form $\Delta_M^u$,
$\Delta_D^u$,
for some $u\in\mathbb R^{n,n}$ with $\|u\|_2=1$ and satisfies
\begin{equation*}
d_{\hi}^{\mathcal{P}}(P(\lambda))=d_{\ker}^{\mathcal{S}}(0,M,D),
\end{equation*}
where $d_{\ker}^{\mathcal{S}}(0,M,D)$ is given by Theorem~\mbox{\rm \ref{main}}. In particular, the structured distance to higher index is
bounded by
\begin{equation}\label{polydistbd}
\sqrt{\lambda_{\min} (M^2+D^2)} \leq d_{\hi}^{\mathcal{P}}(P(\lambda)) \leq \sqrt{2\cdot\lambda_{\min} (M^2+D^2)}.
\end{equation}
\item\label{pp3} The matrix polynomial $P(\lambda)$ is in the closure of polynomials in $\Str$ that are unstable
if and only if $M$ and $D$ have a common kernel or $D$ and $K$ have a common kernel. In particular,
the structured distance to instability $d_{\inst}^{\mathcal{P}}(P(\lambda))$ is attained for a perturbation of the form $\Delta_M^u$,
$\Delta_D^u$,
$\Delta_K^u$ with $\Delta_M^u=0$ or $\Delta_K^u=0$ for some $u\in\mathbb R^{n,n}$ with $\|u\|_2=1$ and satisfies
\begin{equation*}
d_{\inst}^{\mathcal{P}}(P(\lambda))=\min\big\{d_{\ker}^{\mathcal{S}}(0,M,D),\,d_{\ker}^{\mathcal{S}}(0,D,K)\big\}
\end{equation*}
where $d_{\ker}^{\mathcal{S}}(0,M,D)$ and $d_{\ker}^{\mathcal{S}}(0,D,K)$ are given by Theorem~\mbox{\rm \ref{main}}.
Moreover, the distance to instability is bounded by
\begin{equation}\label{polydistbd2}
\sqrt{\alpha} \leq d_{\inst}^{\mathcal{P}}(P(\lambda)) \leq \sqrt{2\cdot\alpha},
\end{equation}
where $\alpha:=\min\big\{{\lambda_{\min} (M^2+D^2)},{\lambda_{\min} (D^2+K^2)}\big\}$.
\end{enumerate}
\end{theorem}
\proof
Parts \eqref{pp1} and \eqref{pp2} immediately follow from Theorem~\ref{distances}. For part \eqref{pp3}, we can apply
Theorem~\ref{thm:24.10.19}, if the matrix polynomial $P(\lambda)$ is singular,  regular of index two, or if it is
regular and has a Jordan block of size two at $\lambda=0$, which is equivalent to the reversal $\lambda^2K+\lambda (D-G)+M$
of $P(\lambda)$ having index two. In all cases the statement follows immediately from \eqref{pp2} applied to $P(\lambda)$ and its reversal,
which has the same structure.
\eproof

\subsection{Removing the coefficient $Q$ in dH systems}\label{sec:noDQ}
Due to the multiplicative structure, in the model representation \eqref{dH} the coefficient $Q$ will present difficulties for
the perturbation analysis and it is an open question whether the definition of a dH system needs this term at all.
In this section we will show that the factor $Q$ can be removed and the system \eqref{dHwoQ}
can be reduced to a system with $Q=I_n$.

Suppose first that $Q$ is invertible. Then the system \eqref{dH} is equivalent to the system
\begin{equation}\label{mm5noQ}
Q^\top  E\dot x  = Q^\top  (J-R)Qx.
\end{equation}
Then setting
$\tilde Q=I$, $\tilde E=Q^\top  E$, $\tilde J=Q^\top  JQ$, $\tilde R=Q^\top  R Q$, we see that for the transformed system
the constraints \eqref{mm5-constr} are satisfied.

If  $Q$ is not invertible then, using the singular value decomposition, there exist orthogonal matrices $U\in \Real^{n,n}$ and
$V \in \Real^{n,n}$ such that
\[
U^\top  QV=\mat{cc} Q_{11} & 0\\ 0 & 0\rix,\quad
U^\top  EV=\mat{cc}E_{11} & E_{12} \\ E_{21} & E_{22} \rix,\quad U^\top ( J-R) U=\mat{cc}L_{11} & L_{12} \\ L_{21} & L_{22} \rix
\]
where in all three block matrices the $(1,1)$ block is square of size $r=\rank( Q)$ and $Q_{11}$ is invertible.
Since $Q^\top  E=E^\top  Q$, we get $Q_1^\top  E_{11}=E_{11}^\top  Q_1$ and $E_{12}=0$, and the transformed system is given by a
reduced dH system
\begin{equation}
\mat{cc}E_{11} & 0 \\ E_{21} & E_{22}\rix \mat{c} \dot x_1 \\ \dot x_2\rix
 = \mat{cc}  L_{11} Q_{11} & 0 \\  L_{21} Q_{11} & 0\rix
\mat{c}  x_1 \\ x_2\rix 
\label{splitsystem}
\end{equation}
where $x_1(t)\in\Real^k$, $x_2(t)\in\Real^{n-k}$.

If the pencil $\lambda E-(J-R)Q$ is regular then also the pencil $\lambda E-Q$ is regular and has index at most one, see \cite{MehMW18}.
Then $E_{22}$ is invertible, and therefore the second block-row of the system reads as
\[
\dot x_2 = E_{22}^{-1}(-E_{21} \dot x_1 + L_{21} Q_1 x_1) ,
\]
$x_1$ does not depend on $x_2$ and the variable $x_2$ can be removed from the system. Then the reduced system
\begin{equation}
 E_{11} \dot x_1   = L_{11} Q_{11} x_1,
\label{reducedsystem}
\end{equation}
satisfies the structured assumption~\eqref{mm5-constr} with $Q_{11}$ being invertible,
so we can apply the previous procedure to obtain a system as in~\eqref{mm5noQ} for the reduced system.

The reduced system \eqref{reducedsystem} may have a different Jordan structure at the eigenvalue zero than system (\ref{dH}).
Indeed, dH pencils with singular $Q$ may have Jordan blocks of size two at the eigenvalue zero
(see \cite{MehMW18} for examples), while it follows from Theorem~\ref{thm:singind} that the eigenvalue zero is semisimple if $Q$ is
invertible. A Jordan block of size $2$ at the eigenvalue $0$ would also mean that the system is unstable.

We highlight that the reduction procedure is not advisable if $E_{22}$ is ill-conditioned. Also, note that
due to the fact that the procedure involves nonorthogonal transformations, the distance to singularity may
change considerably during the process. It is an open problem to characterize the distance to singularity for
a system in the form~\eqref{dH}.

When the solution of the second order system~\eqref{second}, or, equivalently, the solution of the quadratic
eigenvalue problem for $P(\lambda)=\lambda^2 M - \lambda(G-D)+K$ is considered, then the classical approach
is the linearization of the problem. As remarked in the proof of Theorem~\ref{thm:24.10.19}, a particular
linearization of $P(\lambda)$ is of the form $\lambda E (J-R)Q$ with
\begin{equation}\label{E^2}
E=\mat{cc}M&0\\ 0&I\rix,\quad J=\mat{cc}G&-I\\ I&0\rix,\quad R=\mat{cc}D&0\\ 0&0\rix,\quad Q=\mat{cc}I&0\\ 0&K\rix
\end{equation}
which corresponds to a system of the form~\eqref{dH}.
In this case we can remove the matrix $Q$ in a simpler way. Since $Q=Q^\top$ one case use $U=V$, and that $\ker Q=\set0\oplus\ker K$
to reduce the system to the form
\[
Q_1:=U^\top Q U= \matp{I & 0\\ 0 &  K_1},\quad U^\top K U= \matp{ K_1 & 0 \\ 0 & 0},
\]
with some symmetric positive  $K_1\in\Real^{k,k}$. As a result we get a so called \emph{trimmed linearization}, see \cite{ByeMX08},
i.e., a pencil $\lambda E_1 - (J_1-R_1)\in\Real^{n+k,n+k}$, where
\begin{equation}\label{E1^2}
 E_1=\mat{cc}M&0\\ 0& K_1\rix,\quad J_1=\mat{cc}G&-K_2^\top\\  K_2&0\rix,\quad  R_1=\mat{cc}D&0\\ 0&0\rix,
\end{equation}
and $K_2=\matp{ K_1 & 0}\in\Real^{k,n}$. Our aim is now to provide results that allow to compare the distances for the trimmed linearisation
$\lambda E_1-(J_1-R_1)$ with the original distances obtained in  Theorem~\ref{thm:quadpoly}.
In the distances for the reduced system we use the matrices
\[
-J_1^2+E_1^2+R_1^2=  \matp{ M^2+G^\top G+ K_2^\top K_2   +D^2    & -GK_2^\top  \\ K_2G & K_2K_2^\top   }
\]
and
\[
E_1^2+R_1^2=  \matp{ M^2 +D^2    & 0  \\ 0 & K_2K_2^\top   }.
\]
Hence, $\lambda_{\min} (-J_1^2+E_1^2+R_1^2)\leq \lambda_{\min} (K_2K_2^\top)=\lambda_{\min}(K_1^2)$ and if $G=0$ the inequality
becomes an equality.   We get immediately the following result.
\begin{proposition}
For the reduced sytem~\eqref{E1^2} one has the following statements.
\begin{enumerate}[\rm (i)]
\item The structured distance to singularity of the pencil $\lambda E_1-(J_1-R_1)\in\mathcal{L}$ satisfies
\begin{equation}\label{lindist}
	d_{\sing}^{\mathcal{L}}(\lambda E_1- (J_1- R_1))\leq \sqrt{ 2 	\cdot\lambda_{\min} (K_1^2)},
\end{equation}
while for $G=0$ we have
\begin{equation}\label{lindist2}
\lambda_{\min} (K_1^2)\leq d_{\sing}^{\mathcal{L}}(\lambda E_1- (J_1- R_1)).
	\end{equation}
\item The structured distance to higher index of the pencil $\lambda E_1-(J_1-R_1)\in\mathcal{L}$ satisfies
\begin{equation}\label{lindist3}
\sqrt{\beta}\leq	d_{\hi}^{\mathcal{L}}(\lambda E_1- (J_1- R_1))\leq \sqrt{ 2 \cdot	\beta},
\end{equation}
where $\beta=\min\{{\lambda_{\min}(M^2+D^2)},{\lambda_{\min}( K_1^2)}   \}$.
\end{enumerate}
\end{proposition}

\section{Conclusions}
Distance problems in linear differential-algebraic systems with dissipative Hamiltonian
structure have been studied. These include the distance to the nearest singular problem, the distance to the nearest high index
problem, and the distance to instability. The characterization of these distances are open problems for general linear
differential-algebraic systems, while we have shown that for dissipative Hamiltonian systems and related matrix polynomials,
explicit characterizations in terms of common null-spaces of several matrices exist.
\bibliographystyle{plain}

\end{document}